\documentclass[11pt,a4paper,srcltx]{amsart}
\usepackage[francais, english]{babel}
\usepackage{amsfonts}
\usepackage{amssymb}
\usepackage{amsmath}
\usepackage{amsthm}
\usepackage{latexsym}
\usepackage{amsbsy}
\usepackage{alltt}
\usepackage{fontenc}
\usepackage[latin1]{inputenc}

\def\C{\mathbb{C}}
\def\R{\mathbb{R}}
\def\N{\mathbb{N}}

\def\N{\mathbb{N}}

\def\1{\bold{1}}
\def\B{\mathbb{B}}

\theoremstyle{definition}

\theoremstyle{plain}
\newtheorem{def/not}[lem]{Definition/Notations}
\newenvironment{proof3.1}
{\noindent {\it{Proof of theorem 3.1}}}{$\Box$ \linebreak[4]}

\begin{document}

\title[Suites d'interpolation pour $B_{\alpha}^{p}(\B^n)$] {SUR LES SUITES D'INTERPOLATION POUR\\
LES ESPACES DE BERGMAN A POIDS\\
DANS LA BOULE DE $\C^n$}
\footnote{  Ce travail a b\'{e}n\'{e}fici\'{e} du
soutien de l'Action Int\`{e}gr\'{e}e AI 180 MA et du Projet
d'Appui \`{a} la Recherche Scientifique PARS I07. }
\author{ Abdelkader El HASNAOUI}
\address{C.P.R Souissi de Rabat, Avenue Allal Al Fassi, Al Irfane, Rabat
, Maroc   }
\email{e$\_$abdelkader@hotmail.com}
\subjclass[2000]{32A25,
Secondary: 32A36, 65B10.}
\keywords{Suites d'interpolation, espaces de Bergman à poids de la boule 
de $\C^n $, mesures de Carleson.}

\maketitle

\begin{abstract} Let $A$ be a sequence of points of $\B^n$ the unit
ball of $\C^n$. In terms of interpolating vectorial function (or
Amar's function)[1], we give a necessary condition on $A$ to be
interpolating for the weighted Bergman space $B_{\alpha}^{p}(\B^n)$.
In the particular case of Hardy space $H^p(\B^2)$, this condition is
sufficient no optimal. In the main theorem proof, we resolve
Gleason's problem (vectorial form) in $B_{\alpha}^{p}(\B^n)$.
\end{abstract}
\renewcommand{\abstractname}{Résumé}
\begin{abstract} Soit $A$ une suite de points de la boule unité $\B^n$
 de $\C^n$. En termes de fonction vectorielle ( ou fonction d'Amar)
 [1], nous donnons une condition nécessaire pour que $A$ soit
 d'interpolation pour l'espace de Bergman à poids
 $B_{\alpha}^{p}(\B^n)$. Dans le cas particulier de l'espace de Hardy
 $H^p(\B^2)$, cette condition est suffisante non optimale. La preuve
 du théorème principal contient aussi une résolution du
 problème de Gleason (forme vectorielle) dans $B_{\alpha}^{p}(\B^n)$.
\end{abstract}

\section{Introduction}
Depuis les travaux de K. Seip [10], les suites d'interpolation
pour les espaces de Bergman \`{a} poids dans la boule de ${\C}^n$ 
 ont sucsit\'{e} un int\'{e}r\^{e}t particulier. Dans ce contexte, je définis  une extension de la fonction vectorielle interpolante (ou fonction d'Amar)[1] qui est  un vecteur de $( H^{\infty }(\B^n ))^n $, en son analogue dans $( B_{\alpha}^{p}(\B^n)) ^n $,( $ ( H^{\infty }(\B^n )$ étant l'espace des fonctions holomorphes bornées sur $\B^n$). Les preuves des propriétés qui en résultent sont modelées sur celles de [1] et utilisent les notations M.Jevtic' -  X. Massaneda - P.J. Thomas [6].

\subsection{  Notations et rappels.}
  Pour $z$ et $w$ dans ${\C}^n$, on
  \'ecrit: $$z.w=\sum _{j=1}^{n}z_{j}w_{j}, \;<z,w>=z.\bar{w},\;
  \vert z\vert ^{2}=<z,z>
\;\mbox{et}\;\rho (z)=1-\vert z\vert^{2}. $$  $\phi_a$ est
l'automorphisme involutif qui permute l'origine et $a \in
 {\B}^n := \{ \rho > 0  \}$. Les mesures $d\nu_{n}$ et $d\sigma_{n}$ d\'esignent
  celle de Lebesgue et celle d'aire normalis\'ees  respectivement
  sur
  ${\B}^n $ et $\partial {\B}^n$. Pour $\alpha \in ]0,+\infty[$ et
  $p\in ]0,+\infty[$, on note [6]:
  $$B^p_{\alpha}({\B}^n)=\{ f\in L^p({\B}^n,\rho^{\alpha p-1}d\nu_n); f\;\mbox{ holomorphe}
  \}$$ muni de $\Vert f\Vert_{p,\alpha} := \Vert f\Vert
  _{L^{p}(\rho^{\alpha p-1}d\nu_{n})}$. Lorsque $1<p<+\infty$, son dual est
$B^{q}_{\beta}({\B}^n)$ o\`{u} $1/p + 1/q = 1 $ et $\alpha p
=\beta q $. Pour $\alpha =0$, $B^{p}_{0 }({\B}^n)$ d\'{e}signe la
classe de Hardy $H^p({\B}^n)$ munie de
$\Vert f \Vert _{p} =[\sup\;_{0<r<1}%
\displaystyle {\int } _{\displaystyle{\partial }{\B}^n}\vert
f(rw)\vert ^{p} d\sigma _{n}(w)]^{1/p} =\Vert f^{*}\Vert
_{L^p(\sigma_n)} $,  $f^*$ \'etant la limite radiale de $f$. \\
\indent Soit $A=\{a_{k}\}_{k=1}^{\infty} \subset {\B}^n$ une
suite. Les suites $\lambda = \{\lambda _{k}\}_{k=1}^{\infty}
\subset {\C}$  telles que $\Vert \lambda\Vert_{p, n/p+\alpha } :=
\Vert \{\lambda_{k}\rho (a_{k})^{ n/p+\alpha }\}_{k=1}^{\infty}
\Vert _{\ell ^{p} } <+\infty $ forment l'espace $\ell ^{p}_{
n/p+\alpha}$. $A$ est dite d'interpolation-$B_{\alpha
}^{p}({\B}^n)$ (ou $A \in Int(B_{\alpha }^{p}({\B}^n)$) lorsque :
$$( \lambda \in \ell^p_{n/p+\alpha} ) \Rightarrow (\exists f \in
B^{p}_{\alpha}({\B}^n): f(a_{k})=\lambda_{k},\;k=1,2,...).$$
On peut alors choisir $f$ de fa\c{c}on que $\Vert f \Vert _{p,\alpha }\leq C_{A}\Vert \{\lambda
_{k}\}_{k=1}^{\infty} \Vert _{p, n/p+\alpha }$ o\`{u} $C_{A}$ est
la constante d'interpolation de $A$.\newline
 \indent La condition de Carleson:
\begin{eqnarray}
\mbox{inf} \{\prod _{j\neq k} \mid \phi_{a_{j}}(a_{k})\mid,k\;
\mbox {entier} \geq 0\}
>0
\end{eqnarray}
caract\'{e}rise les suites d'interpolation-$H^{p}({\B}^1)$ pour
$p$ dans $]0,+\infty ]$. Pour $n\geq 2$, (1)  (qui est
suffisante)
n' est pas n\'{e}cessaire [5] pour que  $A$ appartienne
 \`a $Int(H^{\infty }({\B}^n))$. En revanche l'existence de  la
fonction vectorielle $B :{\B}^n \longrightarrow {\C}^n$
interpolante pour $A$, introduite par E. Amar dans [1], est
n\'ecessaire pour que $A$ soit dans $Int(H^{\infty}({\B}^n)$. Elle
est aussi suffisante pour que $A \in Int(H^{p}({\B}^2)$ pour tout
$1\leq p < +\infty $ ([1], Th.1.5).  Comme $Int(H^{p}({\B} ^{2}))
\neq Int(H^{q}({\B} ^{2}))$ pour $p\neq q$, l'{e}xistence de $B$
n'est pas n\'{e}cessaire pour que $A$ soit dans $Int(H^{p}({\B}
^{2}))$ lorsque $p$ est fini . Pour qu'il en soit
 ainsi, nous \'{e}tendons la d\'{e}finition de $B$ (dans le
cadre g\'{e}n\'{e}ral des espaces $%
B^{p}_{\alpha }({\B}^{n})$) de fa\c{c}on \`{a} tenir compte de la
variation de $p$ ;  ce qui se formule comme suit:\newline
\newline
\indent {\bf Definition 1.1} {\it Soient $A=\{a_{k}\}_{k=1
}^{\infty}$ une suite de points de ${\B}^n$, $\alpha \geq 0$ et
 $p>0$ deux r\'{e}els. Une application $B :{\B}^n \longrightarrow
  {\C}^n$  {\it est dite interpolante pour $A$ dans $B^{p}_{\alpha }
({\B}^{n})$ lorsqu'il existe $t_{(n,p,\alpha )}=t>0$ et
$C_{(n,p,\alpha )}=C>0$ deux constantes telles que pour tout
entier $k\geq 1$, il existe une matrice $M_{k}(z)$ \`{a}
\'{e}l\'{e}ments dans $B^{p}_{\alpha }({\B}^{n})$ telle que:
\newline \indent$a)\;B(z)=M_{k}(z)\phi _{a_{k}}(z),\;\;
(z\in {\B}^{n})$, avec $\Vert M _{k}\Vert _{p,\alpha } \leq C$;
\newline \indent $b)\;\vert detM_{k}(z)\vert \geq t$ et $\vert
M^{-1}_{k}(z)\vert \leq C$ pour tout $z\in \{ \vert \phi
_{a_{k}}\vert < t\rho (a_{k})^{n(\frac {n}{p} +\alpha )}\}$;
($\Vert M \Vert$ (respectivement $\vert M \vert $) d\'{e}signe
$sup\{\Vert e \Vert $ (respectivement $\vert e \vert $), $e $
 \'{e}l\'{e}ment de $M \}$).}}\newline
\newline
\indent On retrouve évidemment la fonction originale d'Amar en
mettant $\alpha =0$ et $p=\infty$.
\newline Nous sommes maintenant en mesure d'\'{e}noncer le
r\'{e}sultat principal de ce travail :\newline
\newline
\indent {\bf Th\'{e}or\`{e}me 1.1.} {\it Supposons que $2 < p < +\infty
$ et $\alpha  \in \{0\}\cup [1/p,+\infty[$. Alors toute suite
d'interpolation-$B^{p}_{\alpha }({\B}^{n})$ admet une fonction
vectorielle interpolante dans $B^{p/2}_{2\alpha }({\B}^n
)$}.\newline
 \newline
 \indent
La preuve du Th\'{e}or\`{e}me 1.1 se base sur le Lemme d'extension
lin\'eaire (qui va suivre) et le th\'eor\`eme suivant :
\newline
\newline
\indent {\bf Th\'{e}or\`{e}me 1.2.} {\it  Lorsque $2 \leq p <
+\infty$ et
 $\alpha \in \{0\} \cup [1/p,+\infty[$, il existe une constante
  $L_{(n,p,\alpha )}$ telle que : si
$f_{1},...,f_{N}$ sont dans $B^{p}_{\alpha}({\B}^n)$ et $a$ est
dans $\cap _{j=1}^{N} f_{j}^{-1}(0)$, on peut exhiber
$G_{1},G_{2},...,G_{N}$ dans l'espace produit $(B^{p}_{\alpha}
({\B}^n)^n$ tels que:\newline $\ast \;f_{j}(z)=G_{j}(z).\phi
_{a}(z),\;\;\;\;\;(z\in {\B}^{n}, j=1,2,...,N)$,\newline $\ast
\;\Vert \sum _{j=1}^{N}\vert G_{j}\vert ^{2} \Vert _{p/2,2\alpha }
\leq L_{(n,p,\alpha )}\Vert \sum _{j=1}^{N}\vert f_{j} \vert ^{2}
\Vert _{p/2,2\alpha }.$}\newline
\newline
Lorsque $N=1$ et $a=0$, le Th\'eor\`eme 1.2 r\'esout le probl\`eme
de Gleason ([8], Chap.6) pour l'espace  $B^{p}_{\alpha }({\B}
^{n})$.\\ \\
 \indent Nous obtenons la r\'eciproque du Th\'eor\`eme 1.1 dans
 le cas  $n=2$ et $ \alpha =0$:
\newline
\newline
\indent {\bf Th\'{e}or\`{e}me 1.3.}  {\it Supposons que la suite
  $A=\{a_k\}_{k=1}^{\infty} \subset {\B}^2$ {\it v\'{e}rifie : \\ - La mesure $\mu_A := \sum_k
\rho(a_k)^2 \delta_{a_k}$ est de Carleson ($\delta_a$ est la masse
de Dirac en $a$), \\ - $A$ admet une fonction vectorielle
interpolante dans $H^p({\B}^2)$ o\`u $p\in ]3,+\infty[$. \\ Alors
\`{a} toute suite $\{\lambda_k\}_{k=1}^{\infty}$ de $\ell_{2/p}^p$
(resp. $\ell^{\infty }$),
 on peut associer $f$ dans $H^{p/3}({\B}^2)$ (resp. $H^{p/2}({\B}^2))$ telle que
 $f(a_k)= \lambda_k$, $k=1,2,...$}.}
\newline
\newline
Le cas g\'en\'eral ($n\geq 2$
 et $\alpha \geq 0$) n\'ecessiterait la r\'esolution it\'er\'ee de
 l'\'equation-$\overline{\partial}$, partant d'une (0,n)-forme
 (provenant des donn\'ees) et aboutissant \`a une fonction
  dans $B^p_{\alpha}$.
\newline \newline
\indent {\small {\bf Remerciements}. - J'exprime ma profonde  gratitude au
professeur E. Amar. Je remercie les professeurs B. Jennane, X. Massaneda, P.J. Thomas et A.
Zeriahi.\\}

\section{ Lemmes techniques et preuve du Th\'{e}or\`{e}me 1.2}
\vspace*{0.5cm} \indent\indent{\bf Lemme 2.1.} (d'extension
lin\'{e}aire [S.Drury]). {\it Supposons $\alpha \geq 0$ et $p\geq
2$. Soient $A=\{a_{k}\}_{k=1}^{\infty}$
 d'interpolation-$B^{p}_{\alpha }({\B} ^{n})$ et $C_{A}$ sa
 constante d'interpolation.  Alors pour tout entier $N \geq 1$,
 il existe $\beta_{1},\beta {_2 },...,\beta_{N}$ dans $B^{p}_{\alpha
}({\B} ^{n})$ v\'{e}rifiant :\newline $\ast \; \beta _{j}(a_{k})=
\delta _{jk} \;\;\;(symbole\;de\; Kronecker)$,\newline
$\ast\;\Vert \sum_{j=1}^{N} \vert \beta _{j} \vert ^{2} \Vert
_{p/2,2\alpha  } \leq  C_{1} := C_{A}^{2}[\sum _{k=1}^{+\infty
}\rho (a_{k})^{n+\alpha p}]^{2/p}$}.
\newline
\newline
\indent{\bf Preuve}. Ce lemme fut utilis\'{e} dans [1]  avec
$\alpha =0$. On applique  l'\'{e}galit\'{e} de Plancherel  $\Vert
\hat {g_{z}} \Vert _{ L^{2}(G^{\prime}_{N},d\mu'_{N})} =
\Vert g _{z} \Vert _{L^{2}(G_{N},d\mu_{N})}$, aux donn\'{e}es suivantes :%
\newline
$\star\;\;G_{N}=\{1,\lambda,...,\lambda ^{N-1}\}$ est le groupe
multiplicatif
 engendr%
\'{e}  par $\lambda = exp(2i\pi /N)$ et  $d\mu _{N}:=\frac {1}{N}
\sum_{j=1}^{N} \delta _{\lambda ^{j}}$ est sa mesure de Haar ($\delta
_{x} $ est la masse de Dirac au point $x$).\newline
$\star\;\;G^{^{\prime}}_{N}$ est le groupe (cyclique \`{a} N
\'{e}l\'{e}ments) dual de $G_{N}$ de mesure duale $d\mu'_{N} =
\sum_{j=1}^{N}\delta_{\gamma ^{j}}$ o\`{u} $\gamma ^{j}$ est la
$j^{\grave{e}me}_{=}$ puissance du caract\`{e}re $\gamma$
d\'{e}fini par $\gamma (\lambda )=\lambda $.\newline $\star
\;\;g_{z}:G_{N}\longrightarrow \hbox{C\kern -.58em {\raise .54ex
\hbox{$\scriptscriptstyle |$}}
  \kern-.55em {\raise .53ex \hbox{$\scriptscriptstyle |$}} }$ est
  associ\'{e}e \`{a} $z\in {\B}^n$ par $%
g_{z}(\lambda^{j}):=g_{j}(z)$ ($j=1,...,N$) o\`{u} $g_{j}$
r\'{e}sout le probl\`{e}me d'interpolation :
 $\; g_{j}\in
B^{p}_{\alpha }({\B} ^{n})$,
 $\;g_{j}(a_{k})=\lambda ^{jk}$ et
 $\; \Vert g_{j} \Vert _{p,\alpha } \leq
C_{A}(\sum _{k=1}^{\infty }\rho (a_{k}) ^{n+\alpha p})^{1/p}.$\\
Pour $j$ dans $\{1,...,N\}$, les transform\'{e}es de Fourier
$\beta _{j}(z) :=\hat {g_{z}}(\gamma ^{j})$ sont bonnes pour le
lemme.\newline
\newline
\indent {\bf Lemme 2.2.}  {\it  Soient $p\in ]1,+\infty[$ et
$\alpha \in\{0\}\bigcup [1/p,\infty [$. Alors \`{a} $f$ dans
$B^{p}_{\alpha }({\B}^{n})$ et $a$ dans $f^{-1}(0)$, on peut
associer $G_{a}$ dans l'epace $(B^{p}_{\alpha } ({\B}^{n}))^{n}$
de fa\c{c}on que :\newline $ (i)\;\;f(z)=G_{a}(z).\phi _{a}(z),
\;\;\;\;\;(z\in {\B}^{n}),\\ (ii)\;\;\Vert \vert G_{a} \vert \Vert
_{p,\alpha } \leq K_{(n,p,\alpha )} \Vert f \Vert _{p,\alpha }$.
\newline o\`{u} $K_{(n,p,\alpha )}$ est une constante ne
d\'{e}pendant que de $n,p$ et $\alpha$ et $\Vert \vert G_{a} \vert
\Vert _{p,0}:=\Vert \vert G_{a}^{*} \vert \Vert
_{L^p(d\sigma_n)}$.}
\newline
\newline
\indent {\bf Remarque 2.1.} (a)  Pour $a\in{\B}^n$ et
$k_{a}(z):=(\frac {\rho (a)}{(1-<z,a>)^{2}})^{n/p+\alpha }$,
l'op\'{e}rateur  $T _{a}(f):=k_{a} \times (f\circ \phi_{a})$
conserve la norme $\Vert \; \Vert _{p,\alpha }$ de $L^{p}_{\alpha}
({\B}^{n})$ et induit une isom\'{e}trie sur $B^{p}_{\alpha
}({\B}^n)$ ([6], Lemme1.7). Si $G=(g_{1},...,g_{m})$ est un
vecteur de $(L^{p}_{\alpha}({\B}^{n})) ^{m}$,  on a $\vert
T_a(G)\vert = \vert T_a(\vert G \vert)\vert$ o\`u
$T_a(G):=(T_a(g_{1}),...,T_a(g_{m}))$.\\ \indent(b) Lorsque $l$
est entier $\geq1$ et $P_n: {\C}^{n+l}\rightarrow {\C}^n$ est la
projection orthogonale, nous disposons de l'\'egalit\'e de Forelli
([8], p.14) :
 \newline $(F)$\hspace*{2cm}$  \displaystyle{\int
}_{\displaystyle{\partial}{\B}^{n+l}} g\circ P_{n} d\sigma_{n+l}
 =  C^{n+l-1}_{n} \displaystyle{\int }_{{\B}^n}\rho (w)^{l-1}g(w)
 d\nu_{n}(w)$,\\
 \hspace*{2.5cm}$\mbox{pour}\;g\in L^1(\nu_n)\;\mbox{et}\;C^{n+l-1}_{n}:=
 \frac{(n+l-1)!}{n!(l-1)!}$ .\newline
 L'\'egalit\'ee dans $(F)$ demeure vraie dans la situation
 suivante : \newline
$(Fbis)$\hspace*{1.5cm}$  g \;\;\mbox{continue
dans}\;\;{\B}^{n}$,\\ \hspace*{2.6cm}$ g\rho^{l-1}\in
L^{1}(\nu_{n})\; \mbox{et}\; g\circ P_{n}\in
L^{1}(\sigma_{n+l})$.\newline
 ($\chi_{r}(z)$ \'etant $1$ si $|z|<r<1$ et $0$ sinon, on met $g\chi_{r}$
 dans $(F)$ et on fait tendre $r$ vers $1^{-}$).
\newline
\newline
\indent {\bf Preuve du Lemme 2.2.} La Remarque 2.1(a) permet de
supposer $a=0$. En effet, le fait que $T_a(f)(0)=0$ nous donne
$G_{0} \in (B^{p}_{\alpha }{\B}^{n}))^{n}$ avec
$T_a(f)(w)=G_{0}(w).w$
 et $\Vert \vert G_{0}\vert \Vert _{p,\alpha }\leq K_{(n,p,\alpha )}
 \Vert T_a(f) \Vert _{p,\alpha }$. Comme $
k_{a}(\phi _{a}(z)) = (k_{a}(z))^{-1}$ ([8], Th.2.2.2), on prend
$w=\phi _{a}(z)$ pour voir que le vecteur $G_{a}:= T_a(G_{0})$
v\'{e}rifie ($3*$) du Lemme 2.2.\newline \indent $\underline{\mbox
{Cas}\; \alpha =0}$. Puisque $f \in H^{p}({\B}^{n})$ et $f(0)=0$,
la solution d'Ahern-Schneider au probl\`{e}me de Gleason
donn\'{e}e au Chapitre 6 dans [8] est le vecteur $G_{0}$ voulu
dans le Lemme 2.2. Sa composante de rang $k$ est donn\'{e}e par
:\newline $(n)$\hspace*{1.5cm}$ g_{k}^{(n)}(z)= \displaystyle{\int
}_{\displaystyle{\partial } {\B}^{n}} \frac {1}{(1-<z,\xi >)^{n}}
\frac {1-(1-<z,\xi>)^{n}}{<z,\xi >}\bar{\xi _{k} }f(\xi )d\sigma
_{n}(\xi )$\newline  \hspace*{2.cm}avec $f(z) = G_{0}(z).z$ et
$\Vert \vert G_{0}\vert \Vert _{L^{p}(\sigma _{n})}\leq
K_{(n,p,0)}\Vert f \Vert
_{\displaystyle{H^{p}}({\B}^{n})}$.\newline \indent
$\underline{\mbox{ Cas}\; \alpha =l/p\; \mbox{o\`{u}}\; l
\;\mbox{est entier}\; \geq 1 }$. Fixons $f$ dans
$B^{p}_{l/p}({\B}^n)$ telle que $f(0)=0$. Puisque $f\circ P_{n}\in
H^p({\B}^{n+l})$ ([6], Lemme1.2) et $f\circ P_{n}(0)=0$,
l'analogue $(n+l)$ de $(n)$ donne :
\begin{eqnarray}
f \circ P_{n} (z,w)& = &\tilde{G}_{0}(z,w).(z,w),\;\;(z,w) \in
{\B}^{n+l},\;z\in{\C}^l\\& \mbox{o\`u}& \tilde{G}_{0} \in
(H^p({\B}^{n+l}))^{n+l}. \nonumber
\end{eqnarray}
La composante de rang $k$ du vecteur $\tilde{G}_{0}(z,0)$
s'\'ecrit :
\begin{eqnarray}
\tilde{g}_{k}^{(n+l)}(z,0) &=& \int _{\displaystyle{\partial}%
{\B}^{n+l}} h_{z}\circ P_{n}(w,\xi ) d\sigma _{n+l}(w,\xi
),\;\;\;\;(w\in \hbox{C\kern -.58em {\raise .54ex
\hbox{$\scriptscriptstyle |$}}
  \kern-.55em {\raise .53ex \hbox{$\scriptscriptstyle |$}} }^{n}) \\
\mbox{o\`u}\;\;h_{z}(w)&:=& \frac{1}{(1-<z,w>)^{n+l}}
\frac{1-(1-<z,w>)^{n+l}%
}{<z,w>} \bar{w_{k}}f(w).
\end{eqnarray}
\noindent Pour $z$ fix\'e dans ${\B}^n$, on a $h_{z}(w)\rho
(w)^{l-1}\in L^1(d\nu_n)$ et gr\^ace \`a (Fbis), (F) s'applique:
\begin{eqnarray}
\tilde{g}^{(n+l)}_{k}(z,0) &=& C^{n+l-1}_{n}\int_{{\B}%
^{n}} \frac{\rho (w)^{l-1}}{(1-<z,w>)^{n+l}}\\ &\times &
\frac{1-(1-<z,w>)^{n+l}}{<z,w>} \bar{w_{k}} f(w) d\nu_{n}(w).
\nonumber
\end{eqnarray}

\indent $\underline{\mbox{Cas}\;\alpha \in [1/p,+\infty [ }$. Soit
$f\in B^p_{\alpha}({\B}^n)$ avec $f(0)=0$.   (5) nous sugg\'{e}re
de poser:
\begin{eqnarray}
g_{k}(z) &:=& \gamma (n,p,\alpha ) \int_{{\B}^{n}}
 \frac{%
\rho (w)^{\alpha p-1}}{(1-<z,w>)^{n+\alpha p}}\times \\
 & &\frac{1-(1-<z,w>)^{n+\alpha p}}{<z,w>} \overline{w_{k}} f(w)
d\nu_{n}(w),\nonumber
\end{eqnarray}
o\`u $\gamma (n,p,\alpha ):= \frac{\Gamma (n+\alpha p)}{%
\Gamma (\alpha p)\Gamma (n+1)}.$
 Montrons que le vecteur
\begin{eqnarray}
G(z):= (g_{1}(z),...,g_{n}(z))
\end{eqnarray}
\noindent \noindent v\'{e}rifie effectivement:
\begin{eqnarray}
(i)\;\;\;\;\;\;\;\;\;f(z)&=&G(z).z\;, \nonumber \\ (ii)\;\;\Vert
\vert G\vert \Vert _{p,\alpha} &\leq& K_{(n,p,\alpha )}\Vert f
\Vert _{p,\alpha }\;.\nonumber
\end{eqnarray}
On a gr\^ace à (6):
\begin{eqnarray}
G(z).z &=& \sum_{k=1}^{N}g_{k}(z)z_{k}  \nonumber \\
&=&\displaystyle{\int}_{\hbox{I\kern-.2em\hbox{B}}^{n}}f(w)\overline
{K_{z}(w)}\rho (w)^{\alpha p-1}d\nu _{n}(w) \nonumber \\
&-&\gamma _{(n,p,\alpha )}\displaystyle{\int}_{\hbox{I\kern-.2em\hbox{B}}%
^{n}}f(w)\rho (w)^{\alpha p-1}d\nu _{n}(w) \\
\mbox{o\`u} \;\; K_{z}(w)&:=& \gamma_{(n,p,\alpha )}\frac{1}{%
(1-<w,z>)^{n+\alpha p}}
\end{eqnarray}
\noindent La premi\`{e}re int\'{e}grale dans (8) vaut $f(z)$
puisque le noyau $K_{z}(w)$ est reproduisant pour l'espace
$B^{p}_{\alpha }({\B}^n)$ ([6], p.514); alors que la seconde est
nulle puisque $f(0)=0$(l'\'{e}crire en coordonn\'{e}es polaires);
d'o\`{u} le point $(i)$.\newline\indent Montrons $(ii)$.
L'\'{e}galit\'{e} (6) s'\'{e}crit aussi:
\begin{eqnarray} g_{k}(z)\rho
(z)^{\alpha -1/p}& =& \gamma _{(n,p,\alpha )}
\displaystyle{\int}_{{\B}^{n}}\frac{\rho(w)^{s}}{(1-<z,w>)^{n+1+s}}
\lambda _{k}(z,w)f(w)\times \nonumber \\ & & \rho (w)^{\alpha -1/p} d\nu_{n}(w)
\nonumber
\end{eqnarray}
o\`{u} \hspace*{1cm}$ \lambda_{k}(z,w):=(\frac{\rho
(z)}{1-<z,w>})^{\alpha -1/p} \;\frac{1-(1-<z,w>)^{n+\alpha
p}}{<z,w>}\overline{w_{k}}$  \\
 et  \hspace*{2cm} $\;s:= \alpha (p-1)+ 1/p -1 .$ \newline
Puisque $\frac{1-(1-<z,w>)^{n+\alpha p}}{<z,w>} = -(n+\alpha p)
+O(\vert <z,w> \vert)$ et $\frac{\rho (z)} {\vert 1-<z,w>
\vert}\leq 2$, $\lambda _{k} $ est born\'{e}e sur ${\B}^{n}\times
{\B}^{n}$. Le Th\'{e}or\`{e}me de Forelli-Rudin ([8],Chap.7)
implique alors que $\Vert g_{k} \Vert _{p,\alpha } \leq
K_{(n,p,\alpha)}\Vert f \Vert _{p,\alpha }$ (puisque $p\geq 1$);
ce qui montre $(ii)$.\newline\newline \indent {\bf Preuve du
Th\'{e}or\`{e}me 1.2}.  Elle consiste \`{a}
 prouver le lemme suivant :\newline
\newline
\indent {\bf Lemme 2.3}. {\it Soient $\{ f_{j}\}_{j=1}^{N}$ une
partie finie de $B_{\alpha }^{p}({\B}^{n})$ et $a\in \cap
_{j=1}^{N} f_{j}^{-1}(0)$. Les vecteurs $G_{j}$ de $(B^{p}_{\alpha
} ({\B}^{n}))^n$ associ\'{e}s respectivement  aux fonctions
$f_{j}$ par le Lemme} 2.2 {\it v\'{e}rifient pour $2\leq
p<+\infty$ et $\alpha\in \{0\}\bigcup [1/p,+\infty[$} :
\begin{eqnarray}
\Vert \sum_{j=1}^{N} \vert G_{j} \vert ^{2} \Vert _{p/2,2\alpha }
\leq L_{(n,p,\alpha )}\Vert \sum_{j=1}^{N} \vert f_{j} \vert ^{2}
\Vert _{p/2,2\alpha}
\end{eqnarray}
\newline
\indent  Les dualit\'{e}s suivantes sont valables
 pour $p\in ]1,+\infty [$ et $\frac{1}{p}+\frac{1}{q}=1$ :
\newline \indent ${\bf d}_{1}$) Pour $\alpha $ et $\beta $ dans
$]0,+\infty[$ tels que $ \alpha p=\beta q$,  le produit
\[
<U,V>_{N}\; :=
\sum_{j=1}^{N}\int_{{\B}^n}u_j\overline{v_j}\rho^{\alpha
p-1}d\nu_n
\]
 met en dualit\'{e} les deux espaces $(L^{p}_{\alpha }(d\nu_{n}))^{N}$
 et $((L^{q}_{\beta }(d\nu_{n}))^{N}$ munis respectivement des normes 
$\Vert U \Vert _{p,\alpha ,N} :=
\Vert (\sum_{j=1}^{N} \vert u_{j} \vert ^{2} )^{1/2} \Vert
_{p,\alpha }$ et $\Vert V \Vert _{q,\beta ,N}$.\newline \indent
${\bf d}_{2}$) Pour $\alpha =0 $, le produit
\[
<W,H>^{\ast }_{N}\; := \sum_{1}^{N}
\int_{\displaystyle{\partial}{\B}^{n}}
w_j^{\ast}\overline{h^{\ast}_j} d\sigma _{n}
\]
met en dualit\'{e} $(H^{p}({\B}^{n}))^{N}$ et
$(H^{q}({\B}^{n}))^{N}$ munis respectivement des normes
$\Vert W \Vert _{p,N} := \Vert  \vert W \vert  \Vert
_{L^{p}(\sigma _{n})}$ et $\Vert H \Vert _{q,N}$.
 \newline \newline \indent {\bf Preuve du
lemme 2.3.} $C_{(n,p,\alpha)}$ d\'esignera toute constante ne
d\'ependant que de $n$, $p$ et $\alpha$.
\newline
 \indent $\underline{Cas\; \alpha \geq 1/p}$. Soit $G_{j}$ le vecteur
 associ\'{e} \`{a} la donn\'{e}e $f_{j}$
par le Lemme 2.2. Sa kième composante s'écrit grâce à (6):
\begin{eqnarray}
g^{k}_{j}(z)&=&\gamma_{(n,p,\alpha )} \int_{{\B}^{n}} \frac{\rho
(w)^{\alpha p-1}}{(1-<z,w>)^{n+\alpha p}} B_{k}(z,w)f_{j}
(w)d%
\nu_{n}(w)\\ \mbox{o\`{u}}\;\vert B_{k} \vert &\leq
&C_{(n,p,\alpha)}.\nonumber
\end{eqnarray}
Fixons $k$ dans $\{1,...,n\}$ et notons $F:=(f_{1},...,f_{N})$,
$G:=(g_{1},...,g_{N})$ o\`u $g_{j}:=g^{k}_{j}$. Soit
$H=(h_{1},...,h_{N})$ quelconque dans le dual $(L^{q}_{\beta }
({\B}^{n}))^{N}$. Gr\^{a}ce \`{a} (11) on a :
\begin{eqnarray*}
<G,H>_{N}&=& \int_{{\B}^{n}} \sum_{j=1}^{N}g_{j}
(z)\overline{h_{j}(z)} \rho (z)^{\alpha p-1}d\nu_{n}(z)  \nonumber
\\ &=& \gamma _{(n,p,\alpha
)}\int_{{\B}^{n}}[\int_{{\B}^{n}}\frac{\rho (w)^{\alpha
p-1}}{(1-<z,w>)^{n+\alpha p}}B_{k}(z,w)\times \nonumber \\ & &
(\sum_{j=1}^{N}f_{j}(w)\overline{h_{j}(z)})d\nu_{n}(w) \;] \rho
(z)^{\alpha p-1}d\nu_{n}(z).
\end{eqnarray*}
L'in\'{e}galit\'{e} de Schwarz et (11) donnent :
\begin{eqnarray}
\vert <G,H>_{N} \vert &\leq & C_{(n,p,\alpha )}\int_{%
{\B}^{n}} [\int_{{\B}^{n}}\frac{%
\rho (w)^{\alpha p-1}}{\vert 1-<z,w> \vert ^{n+\alpha p}}\\
 & \times &\vert F(w)\vert d\nu_{n}(w)]\vert H(z) \vert\rho
(z)^{\alpha p-1}d\nu_{n}(z).\nonumber \\
 &\leq & C_{(n,p,\alpha )} \int_{{\B}^{n}}[\int_{{\B}^{n}}
\frac{\rho (w)^{\alpha p-1-(\alpha -1/p)}} {\vert 1-<z,w> \vert
^{n+\alpha p-(\alpha -1/p)}}  \nonumber \\ &\times &
f(w)t(z,w)d\nu_{n}(w)]h(z)d\nu_{n}(z),
\end{eqnarray}
o\`{u}:
\begin{eqnarray*}
f&:=& \vert F\vert \rho ^{\alpha -1/p}
  \in L^{p}(d\nu_{n})\\
h&:=&\vert H\vert \rho ^{\alpha p-1 -(\alpha -1/p)} \in
L^{q}(d\nu_{n}),\;\mbox{car}\;1/q+1/p= 1 \;\mbox{et}\;\alpha
p=\beta q \\ \;\;\;\;t(z,w)&:=&\vert\frac{\rho (z)}{1-<z,w>}\vert
^{\alpha -\frac{1}{p}} \leq 2^{\alpha -1/p}
\end{eqnarray*}
\noindent La fonction
\[
\eta (z):=\displaystyle{\ \int }_{{\B}^{n}} \frac{\rho
(w)^{s}}{\vert 1-<z,w> \vert ^{n+1+s}}f(w)d\nu_n(w)
\]
est dans $L^{p}(\nu_{n})$ pour $s= \alpha p-1-\alpha +1/p$
 (section 7.1.4 de [8]) et v\'{e}rifie :
   $$\Vert \eta \Vert _{L^{p}(\nu_{n})} \leq C_{(n,p,\alpha
)}\Vert f \Vert _{L^{p}(\nu_{n})}.$$ Ainsi dans (13),
l'in\'{e}galit\'{e} de H\"{o}lder donne :\newline
\begin{eqnarray*}
\vert <G,H>_{N} \vert & \leq & C_{(n,p,\alpha )}\Vert f \Vert
_{L^{p}(\nu_{n})}\Vert h \Vert _{L^{q}(\nu_{n})} \\ & \leq &
---\;\Vert F \Vert _{p,\alpha ,N} \Vert H \Vert _{q,\beta ,N}.
\end{eqnarray*}
La dualit\'e $\bf{d}_1$) permet alors d'avoir (10).
\newline
\indent $\underline{Cas \;\alpha =0}$. Notons
$F=(f_{1},...,f_{N})\in (H^p({\B}^n))^N$ avec $F(0)=(0,...,0)$ et 
$C(z,w)$ le noyau de Cauchy. Pour $j\in\{1,...,N\}$, le vecteur 
$G_{j}$ associ\'{e} \`{a}
$f_{j}$ par le  Lemme 2.2 est donn\'{e}  par sa kième composante
à l'équation (n) (preuve du lemme 2.2):
\begin{eqnarray}
g_{j}^{k}(z) &=& \int_{\displaystyle{\partial}{\B}^{n}}
C(z,w)(\sum_{l=0}^{n-1}<z,w>^l)\overline{w_k}f_{j}(w)d\sigma
_{n}(w) \nonumber \\ &=& \sum_{\theta \in {\N}^n}^{|\theta |\leq
n-1}l_{\theta}z^{\theta}W_{(\theta,k)}(f_j)(z)
\end{eqnarray}
 o\`{u}
\begin{eqnarray*}
 W_{(\theta,k)}(f)(z)&:=& \displaystyle{\int}_{\displaystyle{\partial}%
{\B}^{n}}C(z,\xi ) \overline{{\xi }^{\theta +\theta _{k}}}f(\xi
)d\sigma _{n}(\xi ) \\ \theta ^{k}&=&(0,...,0,1,0,...0),\;
(1\;\mbox{au
 rang}\; k)\\
l_{\theta }&\mbox{est un entier naturel. }&
 \end{eqnarray*}
Pour avoir (10), il suffit donc
 d'\'etablir pour chaques $0\leq\mid\theta\mid\leq n-1$ et $1\leq k\leq n$
  que :
\begin{eqnarray}
\Vert \sum_{j=1}^{N}\vert W_{(\theta,k)}(f_{j})\vert ^{2} \Vert
_{L^{p/2}(\sigma_n)} \leq \Vert \sum_{j=1}^{N}\vert f_{j} \vert
^{2}\; \Vert _{L^{p/2}(\sigma_n)}.
\end{eqnarray}
 Notons
$W=(w_{1},...,w_{N})$ o\`u  $w_{j} = W_{(\theta ,k)}(f_{j})$. Soit
$H=(h_{1},...,h_{N})$ quelconque dans le dual $(H^{q}({\B}
^{n}))^{N}$. On a :
\[
<W,H>^{\ast }_{N} = \sum_{j=1}^{N}\int_{\displaystyle {\partial }%
{\B} ^{n}} w_{j}(z)\overline {h_{j}(z)}d\sigma _{n}
\]
Par d\'{e}finition de $w_{j}$ on a :
\begin{eqnarray*}
\int_{\displaystyle{\partial}{\B}^{n}}w_{j}(z)%
\overline{h_{j}(z)}d\sigma _{n}(z) = \int_{\displaystyle{\partial} %
{\B}^{n}_{z}}[\int_{\displaystyle{\partial} %
{\B} ^{n}_{\xi }} C(z,\xi ) \overline { \xi ^{(\theta + \theta
_{k})} }f_{j}(\xi )d\sigma _{n}(\xi )] \overline {h_{j}(z)}d\sigma
_{n}(z)
\end{eqnarray*}
 \begin{eqnarray*}
&=&\int_{\displaystyle{\partial}{\B} ^{n}_{\xi }}[\overline
{\int_{\displaystyle{\partial} {\B} ^{n}_{z}}C(\xi ,z)
h_{j}(z)d\sigma _{n}(z) \xi ^{\theta +\theta_{k}}}]
f_{j}(\xi)d\sigma _{n}(\xi )  \\ & = &
\int_{\displaystyle{\partial} {\B} ^{n}}f_{j}(\xi )\overline
{h_{j}(\xi )\xi ^{\theta + \theta _{k}}}d\sigma _{n}(\xi ), \\
\end{eqnarray*}
(le noyau $C(\xi ,z)$ de Cauchy est reproduisant pour
$H^{q}({\B}^{n})$). Ainsi :
\begin{eqnarray*}
\vert <W,H> ^{\ast }_{N} \vert & \leq & \int_{\displaystyle{\partial} %
{\B} ^{n}}(\sum_{j=1}^{N}\vert f_{j} \vert
^{2})^{1/2}(\sum_{j=1}^{N}\vert h_{j}\vert ^{2})^{1/2}d\sigma _{n}
\\ & \leq & \Vert \vert F \vert \Vert _{L^p(\sigma_n)}\; \Vert \vert H
\vert \Vert _{L^q(\sigma_n)}.
\end{eqnarray*}
La dualit\'e $\bf{d}_2)$ permet d'en d'\'eduire que
\newline $$\Vert  \vert W \vert ^{2}\Vert _{L^{p/2}(\sigma_n)}  \leq  \Vert
\vert F \vert ^{2}\Vert
_{L^{p/2}(\sigma_n)}\;.\hspace{0.5cm}{\Box}$$
                    \section{ Conditions n\'{e}cessaires }
 \indent{\bf Th\'{e}or\`{e}me 3.1.} {\it Supposons
$2< p < +\infty $ et $\alpha \in\{0\}\bigcup [1/p,+\infty[$. Alors
toute suite d'interpolation-$B_{\alpha}^{p} ({\B}^n)$ admet une
fonction vectorielle interpolante dans
$B^{p/2}_{2\alpha}({\B}^n)$}.\newline
\newline
\indent {\bf Remarque 3.1.} Le noyau $K_z(w)$ (9) reproduisant
pour $B^p_{\alpha}({\B}^n)$ est une fonction (pour $z$ fix\'e)
dans le dual $B^q_{\beta}({\B}^n)$. Il en r\'esulte que toute
suite $ \{B_N\}_{N=1}^{\infty} \subset B^p_{\alpha}({\B}^n)$
v\'erifiant $ \Vert B_N \Vert _{p,\alpha} \leq R$ admet une
sous-suite qui converge (faiblement et) simplement vers $B \in
 B^p_{\alpha}({\B}^n)$ avec $\Vert B \Vert _{p,\alpha} \leq R$.
 \newline\newline
\indent {\bf Preuve du Th\'{e}or\`{e}me 3.1.} Fixons
$A=\{a_k\}_{k=1}^{\infty}$ dans $Int(B^p_{\alpha}({\B}^n))$. La
Remarque 3.1 permet de restreindre la preuve \`a une partie finie
$A_N=\{a_1,...,a_N\}$ qui est aussi dans
$Int(B^p_{\alpha}({\B}^n))$. En effet, si l'on suppose que pour
chaque $N\in {\N}$, on dispose d'une fonction vectorielle $B_N$
interpolante pour $A_N$ dans $B^{p/2}_{2\alpha} ({\B}^n)$ avec des
constantes $t_{(n,p/2,2\alpha)}$ et $C_{(n,p/2,2\alpha)}$
ind\'epen-\\
dantes de $N$, la limite simple $B$ d'une sous-suite de
$\{B_N\}_{N=1}^{\infty}$ est clairement interpolante pour $A$ dans
$B^{p/2}_{2\alpha}({\B}^n)$. Fixons $N$ dans ${\N}$. Soient
$\beta_1,...,\beta_N$ les fonctions de $B^p_{\alpha}({\B}^n)$
associ\'ees \`a $A_N$ par le Lemme 2.1. Nous \'etablirons que la
fonction
\begin{eqnarray}
B_N:= \sum_{j=1}^{N}\beta_j^2\phi_{a_j}
\end{eqnarray}
est interpolante pour $A_N$ dans $B^{p/2}_{2\alpha}({\B}^n)$, avec
des constantes ne d\'ependant pas de $N$. On v\'erifira $a)$ et
$b)$ de la D\'efinition 1.1 pour $k=1$.\\ Pour simplifier les
notations, on pose :
 \noindent$ B:= B_N $\\ $ a\;    :=
a_1 ,\;
 \beta\; := \beta_1 ,\;
 \phi\; := \phi_{a_1} := (\varphi_1,...,\varphi_n)$ \\
$\phi_j := \phi_{a_j} :=
(\varphi_1^j,...,\varphi_n^j)\;\;\mbox{pour}
 \;\;2\leq j\leq N. $ \\
 Puisque $a \in \bigcap _{j=2}^N \beta_j^{-1}(0)$, le
 Th\'{e}or\`{e}me 1.2
  nous donne
  $G_2,...,G_N$ dans $(B^p_{\alpha}({\B}^n))^n$ tels que
 \begin{eqnarray}
    \beta_j(z) &=& G_j(z).\phi(z),\;\;\;\;\;\;\;\;(z\in {\B}^n, 2\leq j
    \leq N);\\
 \Vert \sum_{j=2}^{N}\vert G_j\vert ^2\Vert _{p/2,2\alpha} &\leq &
 L_{(n,p,\alpha)}\Vert \sum_{j=2}^N\vert \beta_j \vert ^2 \Vert
 _{p/2,2\alpha}\;; \nonumber \\ &\leq &C_2 := C_{(n,p,\alpha )}C_1,
 \;\; (C_1\;\mbox{ du Lemme 2.1}).
 \end{eqnarray}
 (16) et (17) permettent d'\'{e}crire :
 \begin{eqnarray}
          B = \beta^2 \phi + \sum_{j=2}^N \beta_j (G_j.\phi)\phi_j.
 \end{eqnarray}
 Les notations :
 \begin{eqnarray*}
G_j& :=& (g_1^j,...,g_n^j)\\
  \nu_{lk}& :=&
 \sum_{j=2}^N\beta_j\varphi^j_l g_k^j\\
  Y& :=&
  (\nu_{kl})_{{\scriptstyle 1\le k\le n
                               \atop \scriptstyle 1\le l\le
                               n}}\\
 I_n& :=& \mbox{la matrice unit\'ee,}
 \end{eqnarray*}
 donnent l'\'ecriture matricielle :
 \begin{eqnarray}
 B&=&M\phi\\ M &=& \beta^2 I_n + Y.
 \end{eqnarray}
 Chaque \'el\'ement $m_{lk}$ de $M$ v\'erifie :
 \begin{eqnarray}
  \mid m_{lk}\mid  \leq  \mid\beta\mid ^2 + \Sigma_{j=2}^{N} \mid G_j
 ^2\mid
 \end{eqnarray}
 et donc
 \begin{eqnarray}
 \Vert m_{lk}\Vert_{p/2,2\alpha} \; \leq  \;C_1 +
 C_2,\;\;(\mbox{constantes de}\; (18)).
 \end{eqnarray}
 Ainsi $M$ v\'erifie le point $a)$ de la d\'efinition 1.1. D'autre
 part son d\'eterminant s'\'ecrit :
 \begin{eqnarray}
 detM = (\beta^2)^n + \sum_{k=1}^{n}s_k(\beta^2)^{n-k}
 \end{eqnarray}
 o\`u $s_k$ est la somme  des $ \frac{n!}{k!(n-k)!}$ mineurs
 principaux d'ordre $k$ de $Y$. La d\'efinition de $Y$ implique :
\begin{eqnarray}
|s_k| &\leq& \frac{n!}{k!(n-k)!}(\sum_{j=2}^{N}|G_j|^2)^k, \\
s_k(a)&=& 0,  \nonumber \\  \|s_k\|_{p/2k,2k\alpha}&\leq&
n!C_2^n.\nonumber
\end{eqnarray}
De m\^eme on a :
\begin{eqnarray}
\beta^2(a) &=& 1,\\ \|(\beta^2)^r\|_{p/2r,2r\alpha} &=&
\|\beta^2\|^r_{p/2,2\alpha},\; \;(r\;\mbox{entier}\;\leq
n),\nonumber\\ &\leq& C_1^n.\nonumber
\end{eqnarray}
Pour la suite nous aurons besoin de l'estimation suivante ([6],
Lemme 1.4) :
\begin{eqnarray}
&&\exists K_{(n,p,\alpha)}\;\mbox{tel que}\;:\nonumber\\ && |h(z)
- h(a) | \leq K_{(n,p,\alpha)}\Vert h \Vert_{p,\alpha}\rho
(a)^{-(\frac{n}{p}+\alpha)}|\phi_a(z)|\\&& \mbox{pour tous :}\;a
\in {\B}^n, \;h \in B^p_{\alpha}, \;z \in \{|\phi_a|< 1/2
\}.\nonumber
\end{eqnarray}
Dans (24) on a :
\begin{eqnarray}
|detM| \geq |\beta^2|^n -\sum_{k=1}^{n}|s_k| |\beta^2|^{n-k}.
\end{eqnarray}
Appliquons (27) successivement \`a $|\beta^2|^n$ et $s_k$, pour
avoir gr\^ace \`a (25) et (26) que:
\begin{eqnarray}
|detM(z)| \geq 1-C_3\rho (a)^{-2n(\frac{n}{p}+\alpha
)}|\phi_a(z)|,\;\;(z\in \{|\phi_a|< 1/2\}).
\end{eqnarray}
De (22), on d\'eduit que tout mineur $d_{n-1}$ d'ordre $n-1$ de
$M$ v\'erifie :
\[
|d_{n-1}| \leq (n-1)![\; \sum_{j=1}^N|\beta_j|^2 + \sum_{j=2}^N
|G_j|^2\;]^{n-1}.
\]
Il en r\'esulte que :
\[
\|d_{n-1}\|_{\frac{p}{2(n-1)},2(n-1)\alpha} \leq
(n-1)!(C_1+C_2)^{n-1},\;\;(C_i\;\mbox{de}\;(18)).
\]
Comme $d_{n-1}(a) = 1$ ou $0$, (27) donne :
\begin{eqnarray}
|d_{n-1}(z)| \leq 1 + C_4\rho (a)^{-2n(\frac{n}{p}+\alpha
)}|\phi_a(z) |,\;\;(z\in \{|\phi_a| < 1/2 \}.
\end{eqnarray}
Maintenant les constantes $t_{(n,p/2,2\alpha)}$ et
$C_{(n,p/2,2\alpha)}$ de la D\'efinition 1.1, voulue dans le
Th\'eor\`eme 3.1, s'obtiennent imm\'ediatement de (29) et (30). Ce
qui finit la preuve.

 \section{ Conditions suffisantes }
 {\bf Th\'{e}or\`{e}me 4.1.}  {\it Supposons que la
suite
  $A=\{a_k\}_{k=1}^{\infty} \subset {\B}^2$ {\it v\'{e}rifie : \\ - La mesure $\mu_A := \sum_k
\rho(a_k)^2 \delta_{a_k}$ est de Carleson ($\delta_a$ est la masse
de Dirac en $a$), \\ - $A$ admet une fonction vectorielle
interpolante dans $H^p({\B}^2)$ o\`u $p\in ]3,+\infty[$. \\ Alors
\`{a} toute suite $\{\lambda_k\}_{k=1}^{\infty}$ de $\ell_{2/p}^p$
(resp. $\ell^{\infty }$),
 on peut associer $f$ dans $H^{p/3}({\B}^2)$ (resp. $H^{p/2}({\B}^2))$ telle que
 $f(a_k)= \lambda_k$, $k=1,2,...$}.}\\
 \\
 \indent On \'etablira  ce th\'{e}or\`{e}me  en deux \'{e}tapes.
 Dans l'\'{e}tape I, on met  la suite finie $\{a_k\}_{k=1}^{N}$
 dans une boule $R{\B}^2 := \{Rz, z\in {\B}^2\}$ o\`u $R<1$. On admet
 ensuite la Proposition 4.1 (qui va suivre) pour r\'esoudre le
 probl\`eme d'interpolation-$H^{p/3}$(resp. $H^{p/2}$) relatif \`a
 la suite $\{R^{-1}a_k\}_{k=1}^N$ avec des constantes
 d'interpolation ne d\'ependant ni de $N$ dans ${\N}$ ni de $R$
 voisin de 1. Cette ind\'ependance permettera de passer aux limites
 lorsque $R$ tend vers $1^-$ et $N$ tend vers $\infty$; le
 Th\'eor\`eme 4.1 en r\'esultera. L'\'etape II  consiste \`a
 prouver la Proposition 4.1. On utilisera la r\'egularit\'e de
 $B_R(z):=B(Rz)$ sur la fermeture $\overline{{\B}^2}$ et
  certaines propriet\'es (mises en annexe) des mesures de Carleson
   g\'en\'eralis\'ees [4] pour estimer des solutions de
   $\overline{\partial}$ et  $\overline{\partial_b}$ dans
   ${\B}^2$.\\
 \begin{center}
                                 {\bf Etape I }
 \end{center}
  \indent L'existence d'une fonction vectorielle
 interpolante
   entraine la s\'eparation faible suivante:  \\
   \\ \indent {\bf Corollaire 4.1} {\it  Si $(a_k)_{k\in {\N}}$ admet une
   fonction vectorielle interpolante dans $B^q_{\beta}({\B}^n)$, il
   existe $\eta = \eta_{(n,q,\beta )} $ tel que  les hyperboules $
   T_j=\{ |\phi_{a_j}|< 2\eta \rho(a_j)^{n(\frac{n}{q}+\beta)}\}$ soient
   deux \`a deux disjointes}.\\
   \\ \indent{\bf Preuve }. D'une part $d_k := \mbox{det}M_k \in B^{q/n}_{n\beta}$
   et $\|d_k\|_{q/n,n\beta}\leq n!C^n$; d'autre part $d_k(a_j)=0$
   ($0\neq \phi_{a_k}(a_j)
   \in kerM_k(a_j)$ pour $j\neq
   k$). Le point $b)$ de la D\'efinition 1.1 et (27) permettent de voir
    que $\eta := t[2K_{(n,q/n,n\beta)}n!C^n]^{-1}$ est bon pour le Lemme
    car $|d_k(z)|\geq t$ et  $|d_k(z)|<t/2$ lorsque $z\in T_k\cap T_j$
    .\hspace*{0.5cm}$\Box$\\
     \\ \indent {\bf Proposition 4.1.}  {\it Soit $\{\lambda_{k}\}_{k=1}^{\infty} \in
     \ell^p_{2/p}$ (resp. $\ell^{\infty}$). Les hypoth\`eses du
     Th\'eor\`eme} 4.1 {\it permettent d'avoir une constante $C^{'}$ telle que
      pour tout entier $N \geq 1$, il existe $R_N \in  ]0,1[$ v\'erifiant : \\
      $(i)\;  \bigcup _{j=1}^{N}\overline {T_j}\subset
      R_N{\B}^2$\newline
   $(ii)\;\forall R \in ]R_N,1[, \;\exists f_{(N,R)} \in
      H^{p/3}({\B}^2)\;(resp.\;H^{p/2}({\B}^2))$ avec,\\
      \hspace*{0.8cm}a) $f_{(N,R)}(R^{-1}a_j)=\lambda_j$, $j\in
      \{1,...,N\}$,\\
      \hspace*{0.8cm}b) $\Vert f_{(N,R)} \Vert_{H^{p/3}\;(resp.\; H^{p/2})} \leq C'$.}\\
    \\ \indent {\bf Preuve du Th\'eor\`eme 4.1.} Soit $N\in {\N}$
 fix\'e et $R_N \in ]0,1[$ le r\'eel de la Proposition 4.1
 (suppos\'ee prouv\'ee). Le dual de $E=H^{p/3}({\B}^2)$ est
 $E'=H^{p/(p-3)}({\B}^2)$ et  sa boule est faiblement compacte. Il
 existe donc une suite $\{\tau_k\}_{k=1}^{\infty} \subset ]R_N,1[$
 et $f_N\in H^{p/3}({\B}^2)$ telles que pour $k$ tendant vers
 $+\infty$, $\tau_k$ tend vers $1$  en croissant  et
 $f_{(N,\tau_k)}$ tend vers  $f_N$ pour la topologie faible $\sigma
 (E,E')$ avec $\Vert f_N \Vert_{H^{p/3}} \leq C'$. Pour $j$ fix\'e
 dans $\{1,...,N\}$, la suite de terme g\'en\'eral $$g_k(w):=
 C(w,\tau^{-1}_k a_j),\;(w\in {\B}^2,\;k\in {\N}),$$ o\`u $C(w,\zeta
 ):= (1-w.\overline{\zeta})^{-2}$  est le noyau de Cauchy, converge
 fortement dans $E'$ vers $g(w):=C(w,a_j)$ puisque
 \begin{eqnarray*}
 \vert g_k(w)-g(w)\vert &=& \vert
 (1-w.\frac{\overline{a_j}}{\tau_k})^{-1}-(1-w.\overline{a_j})^{-1}
 \vert \times \\ && \vert
 (1-w.\frac{\overline{a_j}}{\tau_k})^{-1}+(1-w.\overline{a_j})^{-1}
 \vert \nonumber \\ &\leq & 4(\frac{1}{\tau_k}-1)(1-\frac{\vert
 a_j\vert}{\tau_0})^{-4}.\nonumber
 \end{eqnarray*}
 Le produit $ <f,g>^*_1$ (voir {\bf d}$_2)$)  met en dualit\'e
 $E\ni f$ et $E'\ni g$  et Le noyau de Cauchy est reproduisant pour
 $E^{'}$. Donc : $$ f_N(a_j) =  <f_N,g>^{*}_1 = \lim_ {k\rightarrow
 +\infty} < f_{(N,\tau_k)},g_k>^{*}_1\; =\;\lambda_j .$$ Nous
 obtenons ainsi une suite $(f_N)_{N\in {\N}} \subset H^{p/3}({\B}^2)$
 telle que :  \\
 $f_N(a_j)=\lambda_j,\;\;j\in \{1,...,N\}$ et
 $\Vert f_N\Vert_{H^{p/3}}\leq C'$ ind\'ependante de$N$. \\
  De nouveau la faible compacit\'e de la boule montre que : \\
- $ \exists (f_{N_k})_{k\in {\N}}$ extraite
 de  $(f_{N})_{N\in {\N}}$,\\
- $ \exists f\in H^{p/3}({\B}^2)$\
 telle que  $f_{N_k}\longrightarrow f$
 lorsque $k\longrightarrow +\infty )$ faiblement.  \\
En particulier,
\begin{eqnarray*}
  f(a_j) = <f,C(.,a_j)>^{*}_1 &=&
\lim_{k\longrightarrow\infty}<f_{N_k},C(.,a_j)>^{*}_1 \\
       & =& \lambda_j.
\end{eqnarray*}
 On finit la preuve du Th\'eor\`eme 4.1 en mettant $H^{p/2}$ \`a la
place de  $H^{p/3}$.
\\
\begin{center}
{\bf Etape II }
\end{center}
\vspace{0.2cm}
 \indent\indent{\bf Notations et rappels
 sur les mesures de Carleson}.\\
- Si $(g,h) : {\B}^2\longrightarrow{\C}^2$ et $R\in ]0,1[$, on
note : $g_R(z)=g(Rz)$; $(g,h)_R=(g_R,h_R)$;
$(\partial_zg)(w):=\sum_{k=1}^{2} \frac{\partial g}{\partial
z_k}(w)dz_k $ et
$\partial_z(g,h)(w):=(\partial_zg(w),\partial_zh(w))$.\\ - Dans le
Corollaire 4.1, on note : $\phi_j:=\phi_{a_j}$, $r_j:=
\eta\rho(a_j)^{n(\frac{n}{p}+\alpha)}$, $T_j(r):= \{\mid\phi_j\mid
<r \}$ et $\eta$ est pris dans $]0,1/4[$ tel que
$\overline{T_j(2r_j)}\cap \overline{T_k(2r_k)}=\emptyset$ pour
$j\neq k.$\\ -L'\'ecriture $x\preceq y$ signifie $x\leq \gamma y$
o\`u $\gamma$ est une constante et $(x\approx y)\Leftrightarrow
(x\preceq y\;\mbox{et}\;y\preceq x)$.\newline -
$\mid\sum_{I,J}\omega_{I,J}dz_I \wedge
d\overline{z}_J\mid:=\sum_{I,J}\mid\omega_{I,J}\mid$
 et $\varepsilon_E$ est la fonction valant 1 dans $E$ et 0
 ailleurs.\\
 - $W^1({\B}^n):=\{$mesure
 de Borel $\mu$ dans ${\B}^n$ : $ \vert \mu \vert (\{z\in {\B}^n, \vert 1-\xi
  .\overline{z}\vert < t \}) = O(t^n),\;\ \forall t>0,\;\forall \xi\in\partial
  {\B}^n\}$ est l'espace des mesures de Carleson.\\
 - Pour $\alpha\in]0,1[$, $W^{\alpha }({\B}^n)$ est l'espace des
   mesures de Carleson d'ordre $\alpha $ d\'efini en [4].\newline
   - $W^{\alpha}_{(0,1)}({\B}^2): = \{(0,1)$-forme
   $v$ \`a coefficients continus dans ${\B}^2$ tels que:
 $\rho^{-1/2}(\vert v_1 \vert + \vert v_2 \vert  + \vert \overline{\partial}
\rho \wedge  v\vert )d\nu_2 \in W^{\alpha }({\B}^2)\}$.\\
    - $W^{\alpha }_{(0,2)}({\B}^2): = \{(0,2)$-formes
       $v$ \`a coefficient $m \in
       C^{\infty }({\B}^2)$ tel que $\rho^{1/2}md\nu_2
       \in W^{\alpha }({\B}^2) \}$. \newline
    \indent  Dans le lemme suivant, ces classes  sont connect\'ees par les
         op\'erateurs $\overline{\partial}$ et $\overline{\partial}_b$.
         \newline\newline
     \indent{\bf Lemme 4.1.} {\it Supposons que  $\alpha\in ]0,1[$ et
    $q=\frac{1}{1-\alpha}$ ou $\alpha =1$ et $q=+\infty$. On a}
        \newline
   (i) {\it $(\omega\in W^{\alpha}_{(0,2)}({\B}^2))
         \Rightarrow
          (\exists v\in W^{\alpha}_{(0,1)}({\B}^2) : \overline{\partial}v=
          \omega )$ pour tout $\alpha\in ]0,1]$}.\newline
  (ii) {\it $ (v\in W^{\alpha}_{(0,1)}({\B}^2)\;\mbox{et}
          \;
          \overline{\partial}v=0)\Rightarrow (\exists s\in L^q
          (\partial{\B}^2):
          \overline{\partial}_bs = v) $ pour tout $\alpha \in  ]0,1[$}.
          \newline
 (iii) {\it Pour une mesure de Borel positive $\mu $ dans
          ${\B}^n$, il ya \'equivalence entre :\newline
 \hspace*{0.5cm}(a) $[\displaystyle{\int}_{{\B}^n} \vert f \vert^
 \tau d\mu]^{1/\tau } \leq C_\tau \Vert f \Vert _{H^\tau ({\B}^2)}$
 pour un certain $\tau \in ]0,+\infty[$ et toute $f\in H^\tau
 ({\B}^n)$, \\ \hspace*{0.5cm}(b) $\mu \in W^1({\B}^n)$}.\\ (iv) {\it
 $(\gamma \in W^{\alpha }({\B}^n))
  \Leftrightarrow (\exists \mu_1 \in W^1({\B}^n), \exists
           h\in L^q(\vert \mu_1 \vert) : \gamma =h\mu_1)$ pour tout} $\alpha
           \in ]0,1]$.\newline
 (v) {\it Si $\alpha \geq 1/p, \gamma \in W^{\alpha
 }({\B}^2)\;et\;g\in H^p({\B}^2),\;alors\;g\gamma\in W^{\alpha -1/p
 }({\B}^2)$.}\newline\newline \indent {\bf R\'ef\'erences}. (i) est
 dans le Th\'eor\`eme 3.9 de [1],\newline (ii) est le Th\'eor\`eme
 7 dans [4],\newline (iii) r\'esulte du Th\'eor\`eme 3.1 dans
 [6],\newline (iv) est en page 156 de [3].
 \newline
 Montrons (v). (iv)$\Rightarrow (\gamma =h\mu_1\; \mbox{o\`u}\;
 \mu_1 \in W^1({\B}^2) \;\mbox{et}\; h\in
 L^{\frac{1}{1-\alpha}}(d|\mu_1|))$. $(iii)\Rightarrow g\in L^p(d|
 \mu_1|)$. Donc $gh \in  L^r(d| \mu_1|)$ o\`u $1/r=1/p + 1-\alpha$.
 De nouveau (iv) $\Rightarrow (g\gamma =gh\mu_1) \in
 W^{1-1/p}({\B}^2)$.\\
 \newline \indent {\bf Lemme
 4.2.} (i) {\it La mesure $ d\mu (z):= \sum _{j=0}^{\infty}
  r_j^{-4} \rho (a_j)^{-1} \varepsilon _{T_j(2r_j)}(z)d\nu_2(z)$
  est de Carleson.}\newline
 \newline
  \indent {\bf Preuve}. On utilise
  (iii) du Lemme 4.1 pour $\tau =1$. Si $f\in H^1({\B}^2)$, on a :
  $$ \int_{{\B}^2}\vert  f \vert d\mu = \sum_{j=0}^{\infty}
  r_j^{-4}\rho(a_j)^{-1}
  \int_{T_j(2r_j)}\vert f \vert d\nu_2 \;.$$
  Si $\zeta \in T_j(2r_j)$, (27) implique que :
  $$\vert f(\zeta) -f(a_j) \vert \preceq \Vert f \Vert
  _{H^1}$$
   et l'\'ecriture $\vert f(\zeta)\vert \leq \vert f(\zeta)- f(a_j)
   \vert + \vert
  f(a_j)\vert $ donne :
  $$ \int_{{\B}^2}\vert  f(\zeta ) \vert d\mu (\zeta ) \preceq I_1
  +I_2,$$
  o\`u
  \begin{eqnarray*}
  I_1 = \Vert f \Vert_{H^1}\sum_{j=0}^{\infty} r_j^{-4}\rho (a_j)^{-1}
  \nu_2 (T_j(2r_j)), \\
  I_2 = \sum_{j=0}^{\infty} r_j^{-4}\rho (a_j)^{-1}\vert f(a_j)
  \vert  \nu_2 (T_j(2r_j)).
  \end{eqnarray*}
   Puisque
   $\nu_2 (T_j(2r_j))\preceq r_j^4\rho (a_j)^3 \;
   ([8],  \S2.2.7)$ et $\sum_{j=0}^{\infty}\vert f(a_j) \vert\rho (a_j)^2
   \preceq\Vert f \Vert _{H^1}$
    (car $\sum_{j=0}^{\infty }\rho(a_j)^2\delta_{a_j}$ est
   de Carleson), on obtient que $ I_1 +I_2 \preceq \Vert f\Vert
  _{H^1}$.\hspace{0.9cm}$\Box$\newline
  \newline
  \indent {\bf Lemme 4.3. }  {\it A tout $N$ dans ${\N}$, on peut
  associer $R_N$ dans $]0,1[$ tel que :\\ a) $\bigcup_{j=1}^N
  \overline{T_j(2r_j)}\subset R_N{\B}^2$,\\ b) $\forall R\in ]R_N,1],
  \forall j\in  \{1,...,N\}$, on ait :
    $\vert
  \phi_j(R^{-1}w)\vert < 2\vert \phi_j (w)\vert$
  pour tout $ w \in
  \Omega_j:= \{ \frac{r_j}{\sqrt{2}}<\vert \phi_j\vert < r_j \}. $}\newline
  \newline
 \indent {\bf Preuve.}  Fixons $N$ dans ${\N}$. Soit $R'_N \in
 ]0,1[$ tel que le compact $\bigcup_{k=1}^{N} \overline{T_k(2r_k)}$
  soit dans $R'_N {\B}^2$ et fixons $j\in \{1,...,N\}$. La fonction
 $$\gamma_j(R,w):= 2\frac{\rho (a_j)\rho (w)}{\vert
 1-\overline{a_j}. w\vert ^2} -\frac{\rho (a_j)\rho
 (\frac{w}{R})}{\vert 1-\overline{a_j}.\frac{w}{R}\vert  ^2}$$ est
 uniform\'ement continue sur le compact $[R'_N,2] \times
 \overline{\Omega_j}$. Soit $\epsilon_N :=$ inf $\{
 \frac{r_k^2}{4}, k=1,...,N \}$, ($r_k$ est en notations). Il
 existe $R''_j \in ]R'_N,1[$ tel que : $$ (R\in ]R''_j,1]
 )\Rightarrow (\;\gamma_j (R,w)< \gamma_j (1,w)+\epsilon_N, \forall
 w\in \Omega_j \;).$$ Par ailleurs,
 \begin{eqnarray*}
  \gamma_j
 (1,w)&=& 1-\vert \phi_j(w) \vert ^2, \forall w\in
 {\B}^2,\;\;\;(\mbox{propriet\'e}\;\mbox{des} \;\phi_a,\;a\in
 {\B}^2),\\
      &<& 1-\frac {r_j^2}{4}, \forall w\in \Omega_j,\\
      &<& 1-\epsilon_N, \forall w\in \Omega_j.
 \end{eqnarray*}
 Il en r\'esulte que
 \begin{eqnarray*}
 (R\in ]R''_j,1]) &\Rightarrow & (\gamma_j (R,w) <1, \forall w\in
 \Omega_j),\\
                  &\Rightarrow &( \vert \phi_j (R^{-1}w)\vert <
 2\vert \phi_j (w)\vert , \forall w\in \Omega_j).
 \end{eqnarray*}
 $R_N=\sup_{1\leq j  \leq N}R''_j$ v\'erifie le point $b)$ du
 Lemme.\newline
 \newline
  {\bf Construction de la fonction $f_{(N,R)}$ voulue dans la
 Proposition 4.1}\\ \indent On suppose les hypoth\`eses de la
 Proposition 4.1
  et on fixe : $N$ dans
 ${\N}$, $R_N$ du Lemme 4.3, $R$ dans $]R_N,1[$ et $0\leq \chi\leq
 1$ une fonction de classe $C^{\infty}$ sur ${\R}$ telle que :
  $ \chi
 (x)=
 \left \{
 \begin{array}
 [c]{c}1\;\mbox{si} \;\vert x \vert <1/2 \\
 0\;\mbox{si}\;\vert x \vert > 1.\;\;\;
 \end{array}
 \right.$
 \newline La fonction
 $$F_{(N,R)}(z):= \sum_{j=1}^{N}\lambda_jr_j^{-2}\chi (\vert
 \phi_j(Rz)  \vert^2 )$$ est  lisse au voisinage de
 $\overline{{\B}^2}$ et v\'erifie
 $$ F_{(N,R)}(R^{-1}a_j)=
 \lambda_j, \;\; 1\leq j \leq N.$$ Pour $z\in R^{-1}T_j$ on a :
 \begin{eqnarray*}
 \phi_j(Rz)= A_{j}(Rz)B_R(z),
 \end{eqnarray*}
 o\`u  $A_j$ est l'inverse de $M_j$; il en r\'esulte que : $$
 \overline{\partial}<\phi_j(Rz),\phi_j(Rz)>\; = \;<B_R(z),
 \overline{{}^tA_{j}(Rz)}{}^t (\partial_z\phi_{j,R})(z)>,$$ o\`u
 $\overline{{}^tA}$ est l'adjointe de $A$. Un calcul direct donne (pour tout $z\in {\B}^2$) :
  $$ \overline{\partial_z}F_{(N,R)}(z)=
 B_{1,R}(z)\omega ^1 _{(N,R)}(z) + B_{2,R}(z)\omega^2_{(N,R)}(z) $$
 avec $\left \{
 \begin{array}
[c]{c} (B_1,B_2):=B \\
 B_{j,R}(z):=B_j(Rz)
 \end{array}
 \right .$
 et pour $k\in \{1,2\}$,
 \begin{eqnarray}
 \;\;\;\;\;\; \omega^k_{(N,R)}(z)&:=& R\sum_{j=1}^N\lambda_jr_j^{-2}\chi ^{'}
 (r_j^{-2}\vert \phi_j(Rz)\vert ^{2})C^k_j(Rz) ;
 \end{eqnarray}
$ C^k_j(w)$ \'etant la compposante de rang $k$ du vecteur
  ${}^tA^j(w){}^t\overline{\partial_z\phi_j(w)}$.
Les  (0,1)-formes $\omega^k_{(N,R)}$ sont de classe $C^{\infty}$
au voisinage de $\overline{{\B}^2}$ et on a :
\begin{eqnarray}
\;\;\;\;\;\;\overline{\partial}_z(\omega^1_{(N,R)}(z))&=&-B_{2,R}(z)
\omega^3_{(N,R)}(z) \\
\;\;\;\;\;\;\overline{\partial}_z(\omega^2_{(N,R)}(z))&=&B_{1,R}(z)
\omega^3_{(N,R)}(z) \\ \;\;\mbox{o\`u}\;\;
 \omega^3_{(N,R)}&=& R^2\sum_{j=1}^N\lambda_jr_j^{-4}
\chi ^{''}(r_j^{-2}\vert \phi_j(Rz) \vert^2)\times\\ & &
detA_j(Rz)\overline{J_c\phi_j (Rz)}d\overline{z_1}\wedge
d\overline{z_2}\nonumber
\end{eqnarray}
($J_c$ est le jacobien complexe).\newline Le support strict de la
d\'eriv\'ee $ \chi^{(l)}(r_j^{-2}\vert \phi_j(Rz)\vert ^2)$
d'ordre $l\geq 1$ est dans $R^{-1}\Omega_j$, qui est \`a son tour
inclu dans $T_j(2r_j)$ et les $T_k(2r_k) $ sont deux \`a deux
disjoints. Ainsi si l'on pose
\begin{eqnarray}
\lambda (z):=
\sum_{j=1}^{\infty}\lambda_j\varepsilon_{T_j(2r_j)}(z),\;\;
z\in{\B}^2,
\end{eqnarray}
 les formes dans (31) et (34) s'\'ecrivent :
\begin{eqnarray}
\omega_{(N,R)}^k &=& \lambda (z)[
m_{(N,R)}^{1,k}(z)d\overline{z_1} +
m_{(N,R)}^{2,k}(z)d\overline{z_2}],\;k\in\{1,2\}\\
\omega_{(N,R)}^3 &=& \lambda (z)
m_{(N,R)}^{3}(z)d\overline{z_1}\wedge d\overline{z_2}\;.
\end{eqnarray}
 Pour majorer ces formes uniform\'ement par rapport \`a $N$ et $R$,
 on utilise les in\'egalit\'es suivantes dont les constantes  ne
 d\'ependent ni de $j\in\{1,...,N\}$ ni de $R\in]R_N,1]$ ni de
 $z\in R^{-1}\Omega_j$ :\\ $\ast\;  \rho (z)\leq   \rho (Rz)\preceq
 \rho(a_j)$\\ $\ast\;|J_c\phi_j(Rz)|\preceq \rho (a_j)^{-3/2}$
 \\$\ast\;| (\partial_z \phi_j)(Rz) |\preceq \rho (a_j)^{-1}$, (car
 $\phi_j \in (H^{\infty})^2$)\\$\ast\;|A_j(Rz)|\preceq 1$ \\ $\ast
 \;|\chi^{(l)}(r_j^{-2}|\phi_j(Rz)|^2)|\preceq
 \varepsilon_{R^{-1}\Omega_j}(z)\preceq
 \varepsilon_{T_j(2r_j)}(z)$\\$\ast\;r_k\preceq \rho (a_k)^4$,
 $\forall k\in {\N}$ \\ $\ast\;\rho(z)\approx \rho(a_k)$,  $\forall
 z\in T_k(1/2)$, $\forall k\in {\N}$. \\ Ceci permet de voir que les
 fonctions $m^{l,k}$ de (36) v\'erifient pour chaque $z$ de ${\B}^2$
 :
 \begin{eqnarray}
 &&max [ \;\rho (z)^{-1/2}\vert m_{(N,R)}^{l,k}(z)\vert , \;\rho
 (z)^{1/2}\vert m^3_{(N,R)}(z)\vert \;] \preceq m(z),
 \;l,k\in\{1,2\},\\&& \mbox{o\`u} \;\;\;m(z):=
 \sum_{j=1}^{\infty}r_j^{-4}\rho (a_j)^{-1}\varepsilon
 _{T_j(2r_j)}(z),\;z\in{\B}^2.
 \end{eqnarray}
 Par ailleurs, si $\{\lambda_j\}_{j=1}^{\infty} \in \ell^{p}_{2/p}$
 (resp.$\ell^{\infty}$), on a :
 \begin{eqnarray*}
 \Vert \lambda\Vert^p_{L^p(md\nu_2)}&=&\sum_{j=1}^{\infty}\vert
  \lambda_j\vert^p r_j^{-4}\rho (a_j)^{-1}\nu_2 (Tj(2r_j))\\
 &\preceq & \sum_{j=1}^{\infty}
  \vert \lambda_j\vert^p\rho (a_j)^2,
  \;\;(\mbox{car}\;\nu_2 (Tj(2r_j)) \preceq r_j^4\rho (a_j)^3)\\
   &<& \infty \\
   \end{eqnarray*}
  et gr\^ace au Lemme 4.1, $\lambda (z)m(z)d\nu_2(z)$ est dans
  $W^{1-1/p}({\B}^2)$ (resp. $W^1({\B}^2)$). Ceci avec (36), (37) et
  (38) se r\'esume comme suit :
  \newline
 \newline
\indent {\bf Lemme 4.4.} {\it  Lorsque
$\{\lambda_{j}\}_{j=1}^{\infty} \in \ell ^p_{2/p}$ (resp.
$\ell^{\infty}$), on a :\\ (i) $\lambda (z) m(z) d\nu_2(z)$ est
dans $W^{1-1/p}({\B}^2)$ (resp. $W^1({\B}^2))$}.\\ (ii) {\it Les
(0,1)-formes $\omega^k_{(N,R)}$ sont dans
$W^{1-1/p}_{(0,1)}({\B}^2)$ (resp. $W^{1}_{(0,1)}({\B}^2))$ et
v\'erifient\\  $$ \rho ^{-1/2}\vert \omega^k_{(N,R)}\vert\preceq
\vert \lambda \vert m,\;\;(k=1,2)$$}
\\ (iii) {\it La (0,2)-forme $\omega^3_{(N,R)}$ est dans
$W^{1-1/p}_{(0,2)}({\B}^2)$ (resp. $W^{1}_{(0,2)}({\B}^2))$ et
v\'erifie \\$$\rho^{1/2}\vert \omega^3_{(N,R)} \vert \preceq \vert
\lambda\vert m .$$}\\
 \newline
 \newline
\indent Nous sommes \`a pr\'esent en mesure de r\'esoudre le
$\overline{\partial}$  et le $\overline{\partial}_b$ dans ces
classes de Carleson pour avoir la fonction $f_{(N,R)}$ annonc\'ee.
E. Amar a utilis\'e des noyaux $K_j(z,\zeta )$
 (not\'es $A_j(z,\zeta ) + B_j(z,\zeta)$ dans [1]) pour
 r\'esoudre $$\overline{\partial} v= \omega ,\;\;\;\;\;\;\;(3*)$$
 de sorte que lorsque la donn\'ee  $\omega = l(z)d\overline{z_1}
 \wedge d\overline{z_2}$ est dans $W^{\alpha}_{(0,2)}({\B}^2)$, la
 (0,1)-forme $v$ \`a coefficients $$ v_j(z) := \int_{{\B}^2_{\zeta}}
 K_j(z,\zeta )l(\zeta )d\nu_2(\zeta )$$ est alors dans
 $W^{\alpha}_{(0,1)}(\B^2)$. En fait, il a \'etablit la
 propriet\'e ( $P$) (l\'eg\'erement plus forte) suivante, dont jouissent
 ces noyaux (preuve du Th\'eor\`eme 3.9 en [1]) :\\  \\
\noindent{\bf Proposition}. {\it  Pour $V_l(z):= \frac{1}{\rho
(z)^{1/2}}\displaystyle{\int} _{{\B}_{\zeta }^2} [\vert
K_1(z,\zeta ) \vert + \vert K_2(z,\zeta ) \vert ]l(\zeta ) d\nu_2
(\zeta )$, on a: ($P$)\hspace{1cm} $[\;\rho (\zeta
)^{1/2}l(\zeta
 )d\nu_2(\zeta)\;
 \in W^{\alpha}({\B}^2)\;]\;\;\;\Rightarrow \;\;\;[\;V_l(z)d\nu_2 (z) \in
 W^{\alpha}({\B}^2)\;]. $ } \\ \\
 Mise dans (3$*$),  $\omega_{(N,R)}^3$ du Lemme 4.4 donne :
  \begin{eqnarray}
  \exists v_{(N,R)}\in W^{1-1/p}_{(0,1)}\;(\mbox{resp.}\;
  W^{1}_{(0,1)}) :
  \overline{\partial}v_{(N,R)}= \omega_{(N,R)}^3.
  \end{eqnarray}
 Gr\^ace \`a $(P)$ et le Lemme 4.4,  $l(\zeta)=\vert \lambda (\zeta )\vert m(\zeta
 )$ v\'erifie :
\begin{eqnarray}
 \rho (z)^{-1/2}\vert v_{(N,R)}(z)\vert &\preceq &
V_l(z),\;(z\in{\B}^2),\\   V_l(z)d\nu_2(z) &\in &
W^{(1-1/p)}\;(\mbox{resp.}\;
 W^{1}).\nonumber
\end{eqnarray}
(31) et (34) permettent de voir que les formes
\begin{eqnarray}
\mu^1_{(N,R)} := \omega^1_{(N,R)} + B_{2,R}\; v_{(N,R)}\\
\mu^2_{(N,R)} := \omega^2_{(N,R)} - B_{1,R}\; v_{(N,R)}
\end{eqnarray}
sont :\\ a) $\overline{\partial}$-ferm\'ees gr\^ace \`a (32) et
(33),\\ b) \`a coefficients continus sur $\overline{{\B}^2}$,
($v_{(N,R)}$ l'est gr\^ace \`a la Proposition A de l'annexe A),\\
c) dans $W^{(1-2/p)}_{(0,1)}$ (resp. $W^{(1-1/p)}_{(0,1)}$ ) : on
applique le Lemme 4.2 \`a
 $v_{(N,R)} \in W^{(1-1/p)}_{(0,1)}$ (resp. $
W^{1}_{(0,1)}$) et $B_{j,R} \in H^p$, ainsi que le Lemme 4.4 aux
formes $\omega^k_{(N,R)}$.\newline \indent La solution $
s^k_{(N,R)}$ de Skoda
 ([4]) de l'\'equation
 \begin{eqnarray}
 \overline{\partial}_b s^k_{(N,R)} = \mu^k_{(N,R)}
 \end{eqnarray}
  est continue sur $ \partial{\B}^2$ gr\^ace \`a la proposition A de
  l'annexe A et $\mu^k_{(N,R)}$ est continue sur
  $\overline{{\B}^2}$. On peut alors appliquer la Proposition 2.1
  de ([9], p.239) pour avoir une fonction $S^k_{(N,R)}$ sur ${\B}^2$ telle
  que :\\
 (i) $S^k_{(N,R)}$ est continue sur $\overline{{\B}^2}$, \\
 (ii) $\overline{\partial}S^k_{(N,R)} = \mu^k_{(N,R)}$, \\
 (iii) la restriction de  $S^k_{(N,R)}$ \`a $\partial{\B}^2$ est
 $s^k_{(N,R)}$.\newline
 \indent On pose :
 \begin{eqnarray}
 f_{(N,R)}:= F_{(N,R)}-B_{1,R} S^1_{(N,R)}-B_{2,R} S^2_{(N,R)},
 \end{eqnarray}
 et on montre ce qui suit : \\
 \newline
 \indent {\bf Lemme 4.5.} {\it Sous les hypoth\`eses de la
 Proposition} 4.1, {\it il existe une constante $ C'$ telle que
 $\forall N\in{\N} , \forall R\in ]R_N,1[$, la fonction $ f_{(N,R)}$
 v\'erifie : }\\ (i) {\it $f_{(N,R)}(R^{-1}a_j)=\lambda_j, \; j\in
 \{ 1,...,N\}$,}\\ (ii) $ f_{(N,R)}$ {\it est holomorphe},\\ (iii)
 $\parallel
  f_{(N,R)}\parallel_{H^{p/3}\;(resp.\;H^{p/2})} \leq C'$.\newline
 \newline
\indent {\bf Preuve du Lemme 4.5.} (i)
$B_{k,R}(R^{-1}a_j)=B_k(a_j)=0 $ et $F_{(N,R)}(R^{-1}a_j) \\
 =\lambda_j$ pour $j\in\{1,...,N\}$. \\
(ii) $\overline{\partial}f_{(N,R)} = \overline{\partial}
F_{(N,R)}-B_{1,R}\mu^1_{(N,R)}-B_{2,R}\mu^2_{(N,R)}=0$.
\newline (iii) Pour $k\in \{1,2\}$, $C^{'}_k$ d\'esignera toute
constante ne d\'ependant ni de $N\in {\N}$ ni de $R\in]R_N,1[$. La
fonction $f_{(N,R)}$ est continue sur $\overline{{\B}^2}$ et sa
restriction \`a $\partial{\B}^2$ s'\'ecrit :
\begin{eqnarray}
f^*_{(N,R)}= -B^*_{1,R}\;s^1_{(N,R)}-B^*_{2,R}\;s^2_{(N,R)}.
\end{eqnarray}
Puisque $B_k \in H^p$ et $\|B^*_{k,R}\|_{L^p(\sigma_ 2)}\leq
\|B^*_{k}\|_{L^p(\sigma_ 2)}$,
  il suffit  d'avoir :
 \begin{eqnarray}
\Vert s^k_{(N,R)} \Vert _{L^{p/2}(d\sigma_2)
(\mbox{\scriptsize{resp.} } \;L^p(d\sigma_2) } \;\leq \;C'_k,
 \end{eqnarray}
 pour avoir (iii).
Dans (44), $s^k_{(N,R)}$ s'\'ecrit moyennant les noyaux de Skoda
 ([4], p.20) :
 \begin{eqnarray}
 s^k_{(N,R)}(z)& = &\int_{{\B}^2_\zeta }K(z,\zeta )\wedge
 \mu^k_{(N,R)}(\zeta )\;+\\& & \int_{{\B}^2_\zeta }L(z,\zeta )\wedge
 \frac{\mu^k_{(N,R)}(\zeta )\wedge \overline{\partial}(|\zeta |^2)}
 {\rho (\zeta )^{1/2}}.\nonumber
\end{eqnarray}
La structure  d'espace de nature homog\`ene  que porte le bord
$X=\partial{\B}^2$ a permis d'\'ecrire ([4], Th.6) : \newline -
les noyaux $K$ et $L$ sous la forme $\tilde{K_t}$ et $\tilde{L_t}$
que nous notons indiff\'eremment ici par $P_t$,\newline -
$s^k_{(N,R)}$ comme la somme des balay\'ees (Annexe B) $
P^{*}_{l^k_{(N,R)}d\nu_2}(z)$ de mesures $ l^k_{(N,R)}(\zeta
)d\nu_2(\zeta )$ qui sont dans $ W^{(1-2/p)}$ (resp.
$W^{(1-1/p)}$), car $\mu^k_{(N,R)}(\zeta ) \in
W^{(1-2/p)}_{(0,1)}$ (resp. $W^{(1-1/p)}_{(0,1)}$). Il suffit donc
pour avoir (47), d'\'etablir que :
\begin{eqnarray}
\Vert P^{*}_{l^k_{(N,R)}d\nu_2}
\Vert_{L^{p/2}(d\sigma_2)(\mbox{\scriptsize{resp.}}
\;L^p(d\sigma_2)) }\leq C'_k
\end{eqnarray}
De (48) on obtient que
\begin{eqnarray}
| l^k_{(N,R)}(\zeta ) | \;\preceq \;\rho (\zeta )^{-1/2}|
\mu^k_{(N,R)}(\zeta )| ,\;\;\;(\zeta\in {\B}^2).
\end{eqnarray}
Gr\^ace \`a (41), (42), (43) et  le Lemme 4.4, on d\'eduit  de
(50), que  pour $j \neq k \in \{1,2\}$, la fonction $l(\zeta )=
\;| \lambda (\zeta ) |m(\zeta)$ v\'erifie :
\begin{eqnarray}
&&| l^k_{(N,R)}(\zeta )|\;\preceq \;l(\zeta ) + | B_{j,R}(\zeta )|
V_l(\zeta ),\\ &&ld\nu_2\;\mbox{et}\; V_ld\nu_2\;\in W^{(1-1/p)}
 \;(\mbox{resp}.\;W^{1}).\nonumber
\end{eqnarray}
De m\^eme gr\^ace au Lemme 4.2 et au fait que $B_{j,R}\in
 H^p$ on a :
$$ | B_{j,R}|  V_ld\nu_2 \in
W^{(1-2/p)}\;(\mbox{resp}.\;W^{1-1/p}).$$
 Maintenant les noyaux $P_t$ de Skoda et $P^0_t$ de Hardy-Littlwood
 v\'erifient les conditions  $(H1)$ et $(H'1)$ de l'annexe B. Le
 point (i) de la Proposition B (de cette Annexe) s'applique
 aux deux mesures $l^k_{(N,R)}d\nu_2$ et $( l + | B_{j,R}|
 V_l)d\nu_2$ pour donner gr\^ace \`a (51) :
 \begin{eqnarray}
&& \Vert P^{*}_{l^k_{(N,R)}d\nu_2}\Vert_{L^{p/2}(d\sigma_2)}
 \leq \Vert P^{0*}_{ld\nu_2}\Vert
_{L^{p/2}(d\sigma_2)} +
 \Vert P^{0*}_{\mid B_{j,R}\mid
 V_{l}d\nu_2}\Vert _{L^{p/2}(d\sigma_2)},\\
 &&(\mbox{resp. partout}\; L^p\; \mbox{au lieu
 de}\;L^{p/2}).\nonumber
 \end{eqnarray}
Gr\^ace au Lemme 4.1 on a aussi :\\
 $$\left \{\begin{array}
 [c]{c} V_ld\nu_2 = hdw_1\;\mbox{ o\`u}\; w_1\in
W^1\;\mbox{et}\;h\in L^p(d| w_1| )\;\;(\mbox{resp.}\;
L^{\infty}(d\vert w_1\vert)) \\
  \Vert B_{j,R}\Vert_{L^p(d| w_1| )}
  \;\preceq  \; \Vert B_{j,R}\Vert_{H^p}.
\end{array}
\right.$$
  \noindent Ceci avec la Proposition B de l'Annexe B permet  d'\'ecrire:
\begin{eqnarray*}
   \Vert P^{0*}_{| B_{j,R}|
    V_{l}d\nu_2}\Vert _{L^{p/2}(d\sigma_2)\;
    (\mbox{\scriptsize{resp.}}\;L^{p}(d\sigma_2))}
   & \leq&   \Vert P^{0*}_{| B_{j,R}|\;
             | h| d| w_1| }\Vert
             _{L^{p/2}(d\sigma_2)\;
             (\mbox{\scriptsize{resp.}}\;L^{p}(d\sigma_2))}\\
  & \preceq&  \Vert B_{j,R} \;h \Vert _{L^{p/2}(d|
   w_1|)\;(\mbox{\scriptsize{resp.}}\;L^{p}(d\vert w_1\vert)}\\
 & \preceq&  \Vert B_{j,R} \Vert _{L^{p/2}(d|w_1|)} \times \\
 & & \Vert h\Vert _{L^{p}(d\vert w_1\vert)\;(\mbox{\scriptsize{resp.}}\;L^{\infty}(d\vert w_1\vert)} \\
 &\preceq&  \Vert B_{j,R} \Vert_{H^p}\Vert h\Vert _{L^{p}
(d\vert w_1\vert)\;(\mbox{\scriptsize{resp.}}\;L^{\infty}(d\vert w_1\vert)}\\
 &\preceq&  \Vert B_j \Vert_{H^p}\Vert h\Vert _{L^{p}
(d\vert w_1\vert)\;(\mbox{\scriptsize{resp.}}\;L^{\infty}(d\vert w_1\vert)}
\end{eqnarray*}
 Gr\^ace \`a (52), la constante :\\
  \indent $C'_k: = \; \Vert
P^{0*}_{ld\nu_2}\Vert_{L^{p/2}(d\sigma_2)} +
 \Vert B_j \Vert_{H^p}\Vert h \Vert _{L^p(d|
 w_1|)\;(\mbox{resp.}\; L^{\infty}(d\vert w_1\vert))},\;\;\;
 (k\neq j \in \{1,2\}),$\\
 v\'erifie (49), ce qui ach\`eve la preuve du Lemme 4.5. et
 \'etablit la Proposition 4.1.
\begin{center}
{\bf                        Annexe A }
\end{center}
\indent {\bf Proposition A}. (i)  {\it En} (40), {\it la
(0,1)-forme $v_{N,R}$ r\'esolvant  $\overline{\partial}v_{N,R} =
\omega _{N,R}^{3}$
 est \`{a} coefficients continus sur
$\overline{{\B}^{2}}$}.\\ (ii) {\it En} (44),  {\it la solution
$s^{k}_{N,R}$ de
       $\overline{\partial }_b s_{N,R}^{k} =  \mu _{N,R}^{k}$
est continue sur $\partial {\B}^{2}$}.\\ \\ \indent {\bf Preuve.}
(i) Le support L de $\omega := \omega^3_{N,R}$ est compact dans
      ${\B}^2$. Les coefficients de la solution $v:=v_{N,R}$ qui
      est en page 12 de [1], s'\'ecrivent :
\begin{eqnarray}
v_{j}(z) = \int_{0}^{2\pi}
(\int_{{\B}^{3}_{\zeta}}K_j(z,r_ze^{i\theta} ,\zeta )\wedge \omega
(\zeta^{'}))\frac {d\theta}{2\pi},
\end{eqnarray}
 $$ \mbox {avec}\:\:\: \left \{
\begin{array}
[c]{c}\zeta = (\zeta^{'},\zeta_{3}),\;\; \zeta^{'} \in {\B}^{2} \\
 z\in {\B}^{2},
r_{z}:= (1-\vert z \vert ^{2})^{\frac{1}{2}}\\
 j\in \{1,2\}.
\end{array}
\right.$$ Le d\'enominateur $D(z,r_ze^{i\theta},\zeta)$ commun aux
noyaux $K_{j}(z,r_ze^{i\theta},\zeta )$ est de module
\begin{eqnarray}
\vert D(z,r_{z}e^{i\theta },\zeta )\vert= \vert
1-z.\overline{\zeta^{'}}-r_{z}e^{i\theta}\overline{\zeta_{3}}
\vert^{5}.
\end{eqnarray}
On met le compact L dans $\{\zeta^{'}\in{\B}^{2}: \vert
\zeta^{'}\vert < \delta \} $, avec  $0 < \delta< 1$ pour avoir :
\begin{eqnarray*}
\vert D(z,r_{z}e^{i\theta},\zeta^{'},\zeta_{3})\vert \geq
 (1-(1-\vert z\vert ^{2})^{\frac{1}{2}} - \delta \vert z \vert)^5,
\end{eqnarray*}
 pour tous
  $\left \{
\begin{array}
 [c]{c} z\in\overline{{\B}^{2}}\
(\zeta^{'},\zeta_{3}) \in
 {\B}^{3}\;\;\mbox {avec}\;\; \zeta^{'} \in L\\
 \theta \in [0,2\pi ].
 \end{array}
 \right.$\\
La fonction  $\varphi_{\delta}(t) := 1- (1-t ^{2})^{\frac{1}{2}}
-\delta t$ demeure $ > 0$ sur
 $ ]\frac{\delta}{1 + \delta^{2}} , 1]$,
il existe donc $m_\delta <1$ et $\varepsilon_\delta
  >0$ tels que :
\begin {eqnarray}
   \vert D(z,r_z e^{i\theta} ,\zeta^{'} ,\zeta_3)\vert \geq \varepsilon
    _\delta
\end{eqnarray}
\mbox{pour tous} $ \left \{
\begin{array}
[c]{c} z \in  C_{(m_{\delta} ,1)}:=\{ z \in {{\C} ^{2}}: m_\delta
\leq \vert z \vert \leq 1\},\\ \zeta = (\zeta^{'} , \zeta _{3})
\in {\B}^{3},\;\mbox{avec}\; \zeta^{'}\in L,\\ \theta \in
[0,2\pi].
    \end{array}
    \right.$\\
Ainsi dans (53), chaque coefficient
 $f_j(z,\theta ,\zeta
 )$ de $K_j(z,r_ze^{i\theta },\zeta ) \wedge \omega
 (\zeta^{'})$ est une fonction continue sur
 $C_{(m_\delta,1)}\times [0,2\pi ]\times p_2^{-1}(L)$ o\`{u}
 $p_2: \overline {{\B}^{3}}\ni (\zeta_1,\zeta_2,\zeta_3)\longrightarrow \zeta^{'}=
 (\zeta_1,\zeta_2)\in {\C}^{2}$. En outre, il existe $C_\delta >0$
 tel que :
 $$\vert f_j(z,\theta ,\zeta)\vert < C_\delta \;\;  \mbox{sur}\;\;
 C_{(m_\delta ,1)}\times [0,2\pi ]\times p_2^{-1}(L). $$
 Par cons\'{e}quent, vue comme int\`{e}grale d\'{e}pendant du
 param\`{e}tre $z$, la restriction de la fonction $v_j(z)$ de (53) \`{a}
 $C_{(m_\delta ,1)}$ est continue : on fait tendre $\vert z \vert$
 vers $1$ et on applique le th\'{e}or\`{e}me de la convergence
 domin\'{e}e. Le point (i) de la Proposition A en r\'esulte puisque
 $v$ est d\'ej\`a continue dans ${\B}^2$ comme \'el\'ement de
 $W_{(0,1)}({\B}^2)$.\\
\indent On d\'emontre le point (ii) en utilisant  Le
Th\'{e}or\`eme 2.2 de
 [2] :
 $ \mu_{N,R}^k$ est une (0,1)-forme $ \overline
 {\partial}$-ferm\'{e}e \`{a} coefficients dans
 $L^{\infty}(d\nu_2)$ et ${\B}^2$ est un convexe de type 2 de
  ${\C}^2$; $s_{N,R}^k$ est donc de Lipschitz sur
  $\partial{\B}^2$.\hspace{0.5cm}$\Box$
\vspace{1cm}
\begin{center}
{\bf                        Annexe B }
\end{center}
\vspace{0.5cm}
 On se r\'ef\`ere aux six premi\`eres pages de [4] pour rappeler ce qui
 suit :
 $(X,\delta, \sigma)$ est un espace de nature homog\`ene. $P^0_t$ est le noyau de
Hardy-Littlwood sur $X$. $P_t$ est un noyau sur $X$ et $M$ la
fonction maximale associ\'ee . Elle associe \`a  $g \in L^q(\sigma
)$, la fonction
\begin{eqnarray*}
  Mg(x) &:= & sup_{\{(t,y):\delta (x,y)<t\}}\vert P_tg(y)\vert \\
\mbox{o\`u}\;\;P_tg(y)&:=&\displaystyle {\int}_X
P_t(y,z)g(z)d\sigma (z).
\end{eqnarray*}
On suppose que  $P_t$ v\'erifie les deux conditions :
\begin{eqnarray*}
(H1):&&\forall q \in ]1, \infty [, \exists C_q > 0: \Vert Mf \Vert
_q \leq C_q \Vert f \Vert _q, \forall f\in L^q(\sigma ),\\
 (H'1):&& \mbox{sup}_{(t,x)\in \small{{\R}^+ }\times X}
\displaystyle {\int}_X \vert P_t(x,y) \vert d\sigma (y) < +\infty.
\end{eqnarray*}
\indent Si $w(t,x)$ est  une mesure de Borel born\'ee sur ${\R} ^+
\times X$ de  variation totale $\vert w \vert $, sa balay\'ee par
$P_t$ est la fonction $$P^{*}_{w}(y) := \int _{{\R}^{+}\times X}
P_t(x,y)dw(t,x) $$ et $w\in W^{\alpha}$ signifie que $P^{*}_{|w|}
\in L^{\frac{1}{1-\alpha }}(d\sigma)$, lorsque $\alpha \in
]0,1[$.\\
\\ \indent {\bf Lemme B}. {\it
Supposons que $P_t$  v\'erifie $(H1)$ et $(H'1)$. Soient $\alpha
\in ]0,1[$ et $p=\frac{1}{1-\alpha }\; .$  Il existe alors $C_p> 0
$ :
 $\Vert P^{*}_{w} \Vert _{L^p(\sigma ) }
 \leq  C_p  \Vert P^{0*}_{\vert w \vert} \Vert _{L^{p}(\sigma )}
 ,\;w\in W^{\alpha }.$ }\\ \\
\indent {\bf Preuve.} On note $\parallel\;\; \parallel_p \;:=
\;\parallel\;\;
\parallel_{L^p(\sigma )}$ et  $q = \frac{1}{\alpha }$ le conjugu\'e de
$p$.  Soit $f$  dans  $L^q(\sigma )$.
 Puisque $P^{0*}_{\vert w\vert}\in L^p(\sigma )$, on a d'une part :
 \begin{eqnarray*}
 A &:=&\displaystyle {\int}_X P^{0*}_{\vert w\vert}(y)\vert
 f(y)\vert d\sigma (y)\\
  &\leq & \Vert P^{0*}_{\vert w\vert}\Vert_p \;\Vert f\Vert_{q}\;,
\end{eqnarray*}
d'autre part, la d\'efinition de la balay\'ee et le th\' eor\`eme
de Fubini donnent :
\begin{eqnarray*}
A  \geq \displaystyle {\int}_{{\R}^+\times X}  \vert
 \displaystyle {\int}_X P^0_t(x,y) f(y) d\sigma (y)\;\vert\;d\vert
w\vert (t,x).
\end{eqnarray*}
Ainsi nous avons :
\begin{eqnarray}
 \displaystyle {\int}_{{\R}^+\times X}\vert P^0_t f(x) \vert
 d\vert w \vert (t,x) \;\leq\; \Vert P^{0*}_{\vert w\vert}\Vert_p \;
 \Vert f \Vert _q \;.
 \end{eqnarray}
Puisque $P_t$ v\'erifie $(H1)$, on remplace $f$ par $Mf$ dans (56)
pour avoir :
\begin{eqnarray}
\displaystyle {\int}_{{\R}^+\times X}\vert P^0_t Mf(x) \vert
 d\vert w \vert (t,x)\; \leq \;C_q\Vert P^{0*}_{\vert w\vert}\Vert_p \;
 \Vert f \Vert _q \;.
 \end{eqnarray}
 Comme $\vert P_tf\vert \;\leq \; \vert P^0_t Mf\vert  $ ([4],
 p.6), l'expression explicite de la balay\'ee ([4], p.5) et le
 th\'eor\`eme de Fubini permettent d'en d\'eduire le Lemme B
 puisqu'alors :
\begin{eqnarray*}
 \vert \displaystyle {\int}_X P^*_w(y)f(y) d\sigma (y)\;\vert\;
\leq \;C_q\Vert P^{0*}_{\vert w\vert}\Vert_p \;
 \Vert f \Vert _q, \;\;f \in L^q(\sigma) .
\end{eqnarray*}
 \indent {\bf Proposition B.}  {\it Soient $w_1 \in W^1$ une
mesure positive et }$1< p< +\infty$ .\\ (i) {\it Si $P_t$
v\'erifie $(H1)$  et $(H'1)$ et si $g$ et $h$ sont dans
$L^p(d\sigma )$ avec $|g| \leq  |h|$, alors   }  $\Vert P^{*}_{ g
dw_1}\Vert _{L^p(d\sigma ) }\;\leq \; C_p \Vert P^{0*}_{\vert
h\vert dw_1}\Vert _{L^p(d\sigma )}\;.$\\ (ii) { \it L'op\'erateur
$O:L^p(dw_1)\longrightarrow L^p(d\sigma ) $ d\'efini par $O(h)=
P^{0*}_{hdw_1}$ est born\'e, c-a-d \\$\exists\; C > 0 \;:\;\Vert
P^{0*}_{hdw_1}  \Vert _{L^p(d\sigma )}\; \leq \; C \Vert h \Vert
_{L^p(dw_1)},\; h\in L^p(dw_1).$}\\ \\($O$ est bien d\'efini
puisque $h\in L^p(dw_1) \Longleftrightarrow hdw_1 \in
W^{1-\frac{1}{p}}$ ([3], p.156).\\
\\
\indent {\bf Preuve.} Le Lemme B donne (i) puisqu'on a
 $ P^{0*}_{|g|dw_1} \leq P^{0*}_{|h|dw_1} $. Pour avoir (ii), on
 consid\`ere $g$ dans le dual $L^q(d\sigma )$. On a
 alors:
\begin{eqnarray*}
 \vert \displaystyle {\int}_X P^{0*}_{hdw_1}(y)g(y)d\sigma
(y)\;\vert  &=& \vert \displaystyle {\int}_X[\displaystyle
{\int}_{{\R}^+\times X} P^0_t(x,y)h(t,x)dw_1(t,x)\;]\times \\
& &g(y) d\sigma (y)\vert \\
& \leq & [{\int}_{{\R}^+\times X}\vert P^0_t
g(x)\vert ^{q} dw_1 (t,x)\:]^{\frac{1}{q}}\times\\
&&[{\int}_{{\R}^+\times X}\vert h(t,x)\vert ^{p} dw_1
(t,x)\:]^{\frac{1}{p}}.
\end{eqnarray*}
Gr\^{a}ce  \`{a} une estimation dans ([2], p.8), la premi\`{e}re
des deux derni\`{e}res int\`{e}grales est major\'{e}e par $C\Vert
g \Vert _{L^q(d\sigma )}$ o\`{u} $C$ ne d\'{e}pend pas de
$g$.\hspace*{0.5cm}$\Box$

\end{document}